\documentclass{article}\begin{document}
\centerline{\bf\large On the Complete Evaluation of Jacobi Theta Functions}
\[
\]
\textbf{\centerline{Nikos D. Bagis}}
\textbf{\centerline{Aristotele University of Thessaloniki-AUTH}}
\textbf{\centerline{Address:}}
\textbf{\centerline{Stenimahou 5 Edessa}}
\textbf{\centerline{Pellas 58200, Greece}}
\textbf{\centerline{email-nikosbagis@hotmail.gr}}
\[
\]
\begin{abstract}
Using numerical, theoretical and general methods, we construct evaluation formulas for the Jacobi $\theta$ functions. Some of our results are conjectures, but are verified numerically. 
\end{abstract}

\section{Introduction}

Let $K(x)$ be the complete elliptic integral of the first kind
\begin{equation}
K(x)=\frac{\pi}{2}{}_2F_1\left(\frac{1}{2},\frac{1}{2};1;x^2\right)\textrm{, }|x|<1
\end{equation}
and $k_r$, the elliptic singular modulus, solution of the equation
\begin{equation}
\frac{K\left(\sqrt{1-k_r^2}\right)}{K(k_r)}=\sqrt{r}.
\end{equation} 
When $r$ is positive rational the function $k_r$ take algebraic values.\\
The 3rd and 4th-Jacobian theta function are given by  
\begin{equation}
\vartheta_3(z,q)=\sum^{\infty}_{n=-\infty}q^{n^2}e^{2niz}
\end{equation}
and
\begin{equation}
\vartheta_4(z,q)=\sum^{\infty}_{n=-\infty}(-1)^nq^{n^2}e^{2niz}
\end{equation}
where $|q|<1$ and $z\in \textbf{C}$. Avoiding the above definitions we will use the equivalent notations
\begin{equation}
\theta_3(a,b;q):=\sum^{\infty}_{n=-\infty}q^{an^2+bn}
\end{equation}
and
\begin{equation}
\theta_4(a,b;q):=\sum^{\infty}_{n=-\infty}(-1)^nq^{an^2+bn}
\end{equation}
Also we shall restrict to the case of $a$ positive rational and $b$ general rational.\\ 
In [8] we have shown that if $|q|<1$ and 
\begin{equation}
A(a,p;q):=q^{\frac{p}{12}-\frac{a}{2}+\frac{a^2}{2p}}\prod^{\infty}_{n=0}\left(1-q^{np+a}\right)\left(1-q^{np+p-a}\right),
\end{equation}
with $p>0$, then
\begin{equation}
A(a,p;q)=\frac{q^{\frac{p}{12}-\frac{a}{2}+\frac{a^2}{2p}}}{\eta(q^p)}\theta_4\left(\frac{p}{2},\frac{p}{2}-a;q\right),
\end{equation}
where $\eta$ is the Dedekind-Ramanujan eta function
\begin{equation}
\eta(q):=\prod^{\infty}_{n=1}(1-q^n)\textrm{, }|q|<1. 
\end{equation}
If $q=e^{-\pi\sqrt{r}}$, $r>0$, then 
$$
\eta(q)=2^{1/3}\pi^{-1/2}q^{-1/24}k^{1/12}k^{{*}{1/3}}K(k)^{1/2},\eqno{(9.1)}
$$
where $k=k_r$ and $k^{*}=\sqrt{1-k^2}$.\\
The above identity (8) is simple consequence of the Jacobi triple product formula (see [16]).   

\section{Algebricity}

Having define what we need we state our\\
\\
\textbf{Conjecture.}\\
Let $q=e^{-\pi\sqrt{r}}$ with $r>0$, then for $a,p$ rationals and $p>0$ there always exist algebraic function $Q(x)=Q_{\{a,p\}}(x)$, $0<x<1$, such that 
\begin{equation}
A(a,p,q)=Q_{\{a,p\}}(k_r)\textrm{, }\forall r>0.
\end{equation}
\\
\textbf{Corollary.}\\
If $q=e^{-\pi\sqrt{r}}$ and $r,|a|,p\in\bf Q^{*}_{+}\rm$, then 
\begin{equation}
A(a,p;q)=Algebraic\textrm{  }Number.
\end{equation}
\\

Some verifications of the above Conjectrure have given in [8] for the case of theta functions of the form
\begin{equation}
\sum^{\infty}_{n=-\infty}q^{n^2+mn}\textrm{, } m\in\textbf{Z} 
\end{equation}
Moreover it has been shown that:\\
\textbf{i)} if $m=2s$ (even), then
\begin{equation}
\sum^{\infty}_{n=-\infty}q^{n^2+2sn}=q^{-s^2}\sqrt{\frac{2K(k_r)}{\pi}}
\end{equation}
\textbf{ii)} if $m=2s+1$ (odd), then
\begin{equation}
\sum^{\infty}_{n=-\infty}q^{n^2+(2s+1)n}=2^{5/6}q^{-(2s+1)^2/4}\frac{(k_{11}k_{12}k_{21})^{1/6}}{k_{22}^{1/3}}\sqrt{\frac{K(k_{11})}{\pi}}
\end{equation}
where
$k_{11}=k_r$, $k_{12}=\sqrt{1-k_{11}^2}$, $k_{21}=\frac{2-k_{11}^2-2k_{12}}{k_{11}^2}$, $k_{22}=\sqrt{1-k_{21}^2}$ in view of the evaluation formula (see [8]):
\begin{equation}
\eta(q)^8=\frac{2^{8/3}}{\pi^4}q^{-1/3}k_r^{2/3}(k^*_r)^{8/3}K(k_r)^4.
\end{equation}

Examples of the above conjecture can also be found if we consider the function
\begin{equation}
r=k_i(x):=k^{(-1)}(x)=\left(\frac{K\left(\sqrt{1-x^2}\right)}{K(x)}\right)^2,
\end{equation} 
which is the inverse function of the singular modulus $k_r$. Our method consists of inserting the value $r=k_i\left(\frac{m}{n}\right)$, where $0<m<n$, $m,n$ integers into the form $A(a,p;q)$, and get numerically, using the routine $RootApproximant$ of the program Mathematica, a minimal polynomial which is esentialy the value of an algebraic number. This lead us to conclude that beneath (any theta function) exists minimal polynomials with coefficients rational functions of $k_r$ (in all cases if Conjecture holds). Hence, for every pair of fixed numbers $a,p$, we have a unique algebraic function $Q_{\{a,p\}}(x)$.\\A very easy example to see this someone is with $a=1$ and $p=4$. In this case all the values of $A\left(1,4,e^{-\pi\sqrt{k_i(\frac{m}{n})}}\right)^{24}$ are rationals. With a simple algorithm one can see that
\begin{equation}
A\left(1,4,e^{-\pi\sqrt{k_i\left(r\right)}}\right)^{24}=\frac{16(1-r^2)^2}{r^2}.
\end{equation}
Hence 
\begin{equation}
A\left(1,4,e^{-\pi\sqrt{r}}\right)=\sqrt[12]{\frac{4(1-k_r^2)}{k_r}}
\end{equation} 
and from relation (8):\\
\\
\textbf{Theorem 1.}\\
If $q=e^{-\pi\sqrt{r}}$, $r>0$
\begin{equation}
\theta_4(2,1;q)=\sum^{\infty}_{n=-\infty}(-1)^nq^{2n^2+n}=q^{1/24}\eta(q^4)\sqrt[12]{\frac{4(1-k_r^2)}{k_r}}\textrm{, }\forall r>0.   
\end{equation}      

The continuation follows from the validity of (8) and (10) in rationals and from the fact that every real number is a limit of a rational sequence.\\
For to find the $n-$th modular equation of $A(4,1,q)$ we use Theorem 3 below to get
\begin{equation}
\Pi_n(x)=Q_{\{1,4\}}\left(k \left(n^2k_i\left(Q^{(-1)}_{\{1,4\}}(x)\right)\right)\right),
\end{equation}
where $Q_{\{1,4\}}(x)=\sqrt[12]{\frac{4(1-x^2)}{x}}$ and $Q^{(-1)}_{\{1,4\}}(x)=\frac{1}{8}\left(-x^{12}+\sqrt{64+x^{24}}\right)$.\\
\\

Another example is setting $p=2$ and $a=1/2$, where we find
\begin{equation}
A\left(1/2,2,e^{-\pi\sqrt{r}}\right)=\sqrt[24]{\frac{4(1-k_r)^4}{k_r(1+k_r)^2}}.
\end{equation} 
This is the same theta function as (17), (by changing $q\rightarrow q^{1/2}$). 
For avoiding these cases it is useful to know that:
\begin{equation}
\sum^{\infty}_{n=-\infty}(-1)^nq^{an^2+bn}=\sum^{\infty}_{n=-\infty}(-1)^nq^{an^2-bn}\textrm{, }\sum^{\infty}_{n=-\infty}q^{an^2+bn}=\sum^{\infty}_{n=-\infty}q^{an^2-bn}
\end{equation}
and also if $s$ is positive integer, then
$$
\sum^{\infty}_{n=-\infty}(-1)^nq^{an^2+bn}\textrm{, }\sum^{\infty}_{n=-\infty}(-1)^nq^{asn^2+bsn}
$$
and
\begin{equation}
\sum^{\infty}_{n=-\infty}q^{an^2+bn}\textrm{, }\sum^{\infty}_{n=-\infty}q^{asn^2+bsn},
\end{equation}
are equivalent.\\
\\
\textbf{Theorem 2.}\\
For $q=e^{-\pi\sqrt{r}}$, $r>0$, we have
\begin{equation}
\sum^{\infty}_{n=-\infty}(-1)^nq^{2n^2+3n/2}=q^{-11/96}\eta(q^4)\sqrt[48]{\frac{4\left(1-k_r\right)^4\left(2+k_r-2\sqrt{1+k_r}\right)^{12}}{k_r^{13}\left(1+k_r\right)^2}}
\end{equation}

\section{The Algorithm}

In this section we give the algorithm for finding the expression $Q_{\{a,p\}}\left(k_r^2\right)$ (here we assuming, for simplicity reasons, that $A(a,p;q)=Q_{\{a,p\}}(k_r^2)$ and not $Q_{\{a,p\}}(k_r)$). Our method is based on interpolation. We find a minimal polynomial 
\begin{equation}
P(y,x)=\sum^{N}_{n=0}\sum^{M}_{l=0}a_{nl}y^nx^l,
\end{equation}
such that 
\begin{equation}
P\left(A(a,p;q),m_1(q)\right)=0,
\end{equation} 
with $m_1(q)=k_r^2$ and then solve with respect to $A(a,p;m_1(q))$ (if $P(y,x)=0$ is solvable with respect to $y$).\\
\\  
\textbf{The algorithm} (In Mathematica Program) 
$$
Clear[A]
$$
$$
eta[q]:=QPochhammer[q,q]
$$
$$
A[a,p,q]:=q^{p/12-a/2+a^2/(2 p)}eta[q^p]^{-1} Sum[(-1)^n q^{n^2 p/2+(p-2a) n/2},\{n,-100,100\}]
$$
$$
Clear[q,x,y,u,v]
$$
$$
m[q]:=InverseEllipticNomeQ[q]
$$
$$
x=Series[A[a,p,q]^{12},\{q,0,M\}];
$$
$$
y=Series[m[q],\{q,0,M\}];
$$
$$
t=Table[Coefficient[Sum[c[i,j]x^iy^j,\{i,0,s\},\{j,0,s\}],q^{n}]==0,\{n,1,s^2\}];
$$
$$
rr=Table[c[i,j],\{i,0,s\},\{j,0,s\}];
$$
$$
rr1=Table[u^i v^j,\{i,0,s\},\{j,0,s\}];
$$
$$
mm=Normal[Extract[CoefficientArrays[t//Flatten,rr//Flatten],2]];
$$
$$
m0=Normal[mm];
$$
$$
r1=Take[NullSpace[m0],1] . Flatten[rr]
$$
$$
Take[NullSpace[m0],1] . Flatten[rr1]//Factor
$$

\section{Theoretical Results and Directions}

In this section (here the notation is the traditional i.e.  $A(a,p;q)=Q_{\{a,p\}}(k_r)$), we will try to characterize these functions $Q_{\{a,p\}}(x)$. For this, assume that $\Pi_{n}$ is the $n$-th modular equation of $A(a,p;q)$. then 
\begin{equation}
A(a,p;q^n)=\Pi_{n}\left(A(a,p;q)\right).
\end{equation}
Also assume that our conjecture (relation (10)) hold. Then 
$$
Q_{\{a,p\}}\left(k_{n^2r}\right)=\Pi_{n}\left(Q_{\{a,p\}}(k_r)\right).
$$
By using (16), we get
$$
Q_{\{a,p\}}\left(k_{n^2k_i(x)}\right)=\Pi_n\left(Q_{\{a,p\}}(x)\right).
$$
Setting 
\begin{equation}
S_n(x):=k_{n^2k_i(x)},
\end{equation}
we have the next\\
\\ 
\textbf{Theorem 3.}\\
If the $n$-th modular equation of $A(a,p;q)$ is that of (27),  then 
\begin{equation}
k_{n^2k_i(x)}=S_n(x)=Q_{\{a,p\}}{}^{(-1)}\left(\Pi_{n}\left(Q_{\{a,p\}}(x)\right)\right)\textrm{, }n=2,3,4,...
\end{equation}
If one manage to solve equation (29) with respect to $Q_{\{a,p\}}(x)$ for given $a,p$, then  
\begin{equation}
\sum^{\infty}_{n=-\infty}(-1)^nq^{pn^2/2+(p-2a)n/2}=q^{-\frac{p}{12}+\frac{a}{2}-\frac{a^2}{2p}}\eta(q^p)Q_{\{a,p\}}(k_r), \forall r>0
\end{equation}
and $Q_{\{a,p\}}(x)$ will be a root of a polynomial of degree $\nu=\nu(a,p,x)$. In case $a,p$ are integers and $x$ is rational, then $\nu=\nu(a,p)$ and its coefficients will be integers.\\ 
\\

Note that, in case of rational $x\in(0,1)$ and $a,p$ rational with $0<a,p$, then the degree $\nu$ is independent of $x$ and the minimal polynomial of $Q_{\{a,p\}}(x)$ will have integer coefficients.\\
\\
\textbf{Example.}\\
The 2nd degree modular equation of $A(1,4;q)$ is
\begin{equation}
16 u^8+u^{16}v^8-v^{16}=0.
\end{equation} 
If we solve with respect to $v$, we get $v=P_2(u)$, where $v=A(1,4;q^2)$ and $u=A(1,4;q)$. Moreover
\begin{equation}
P_2(w)=\frac{\left(w^{16}+w^4\sqrt{64+w^{24}}\right)^{1/8}}{2^{1/8}}.
\end{equation}
It is $n=2$. Then hold (see [9])
\begin{equation}
k_{4r}=\frac{1-\sqrt{1-k_r^2}}{1+\sqrt{1-k_r^2}}.
\end{equation}
Hence
\begin{equation}
S_2(x)=k_{4k_i(x)}=\frac{1-\sqrt{1-x^2}}{1+\sqrt{1-x^2}}.
\end{equation}
Finally, we get from the relation (29) of Theorem 3:
\begin{equation}
\frac{\sqrt[8]{Q_{\{1,4\}}(x)^{16}+Q_{\{1,4\}}(x)^4\sqrt{Q_{\{1,4\}}(x)^{24}+64} }}{\sqrt[8]{2}}=Q_{\{1,4\}}\left(\frac{1-\sqrt{1-x^2}}{1+\sqrt{1-x^2}}\right),
\end{equation}
which have indeed a solution 
$$
Q_{\{1,4\}}(x)=\sqrt[12]{\frac{4(1-x^2)}{x}}.
$$  
\\
\textbf{Notes.}\\
i) We note here that function $m_1(q)=k_r^2$, exists in program Mathematica and is called ''InverseEllipticNomeQ[q]''.\\
ii) A useful expansion is
\begin{equation}
k_r=\sqrt{m_1(q)}=4q^{1/2}\exp\left(-4\sum^{\infty}_{n=1}q^n\sum_{d|n}\frac{(-1)^{d+n/d}}{d}\right),
\end{equation}        
where $q=e^{-\pi\sqrt{r}}$, $r>0$.\\
iii) Also $m_1(q)$ can defined with Jacobi null theta functions: If
\begin{equation}
\theta_2(q)=\sum^{\infty}_{n=-\infty}q^{(n+1/2)^2}\textrm{ and }\theta_3(q)=\sum^{\infty}_{n=-\infty}q^{n^2}\textrm{, }|q|<1,
\end{equation}
then
\begin{equation}
m_1(q)=\left(\frac{\theta_2(q)}{\theta_3(q)}\right)^4.
\end{equation}
\\
\textbf{Definition 1.}(see [15])\\
For any smooth function $G$, we define $m_G(x)$ to be such that
\begin{equation}
x=\pi\int^{+\infty}_{\sqrt{m_G(x)}}\eta_D\left(it/2\right)^4G\left(R\left(e^{-\pi t}\right)\right)dt,
\end{equation}
where $R(q)$ is the Rogers-Ramanujan continued fraction and $\eta_D(z)$ is Dedekind's eta function. i.e.
\begin{equation}
R(q)=\frac{q^{1/5}}{1+}\frac{q^1}{1+}\frac{q^2}{1+}\frac{q^3}{1+}\ldots\textrm{, }|q|<1
\end{equation}
and
\begin{equation}
\eta_D(z)=q^{1/24}\prod^{\infty}_{n=1}\left(1-q^n\right)\textrm{, }q=e(z):=e^{2\pi i z}\textrm{, }Im(z)>0.
\end{equation}  
\\
\textbf{Theorem 4.}(see [15])\\
If $y(x)$ is a function defined from the integral equation:
\begin{equation}
5\int^{y(x)}_{0}\frac{G(t)}{t\sqrt[6]{t^{-5}-11-t^5}}dt=x,
\end{equation}
then
\begin{equation}
y(x)=R\left(e^{-\pi\sqrt{m_G(x)}}\right).
\end{equation}
\\

Also if we set  
\begin{equation}
m_G(A):=k_i\left(Q_{\{a,p\}}^{(-1)}(A)\right),
\end{equation} 
then 
\begin{equation}
m_G^{(-1)}\left(n^2m_G(A)\right)=\Pi_n(A)=Q^{*}_{n^2}(A)
\end{equation}
and in view of Theorem 7 of [15] we have
\begin{equation}
y\left(\Pi_n(A)\right)=\Omega_n\left(y(A)\right),
\end{equation}
where $v(x)=R\left(e^{-\pi\sqrt{x}}\right)$ and $\Omega_n(x)=v\left(n^2\cdot v^{(-1)}(x)\right)$. Hence
\begin{equation}
m_G^{(-1)}(n)=Q_{\{a,p\}}(k_n)=\Pi_{\sqrt{n}}\left(m_G^{(-1)}(1)\right)\textrm{, }\forall n\in \textbf{R}_{>0}
\end{equation}
and
\begin{equation}
y(x)=R\left(e^{-\pi\sqrt{k_i\left(Q^{(-1)}_{\{a,p\}}(x)\right)}}\right).
\end{equation}
Note that then $\Pi_{n}(A)$ is that of (27).\\
Also if we set 
\begin{equation}
\theta^{*}(r):=\theta(q):=\theta_{\{a,p\}}(q):=q^{p/12-a/2+a^2/(2p)}\frac{\theta_3\left(\frac{p}{2},\frac{p-2a}{2};q\right)}{\eta\left(q^p\right)},
\end{equation}
where $q=e^{-\pi\sqrt{r}}$, then using Poisson summation formula one can show that
\begin{equation}
\theta_{\{a',p'\}}\left(q'\right)=q^{-\left(p/2-a\right)^2/(2p)}\theta_{\{a,p\}}(q),
\end{equation}
where
\begin{equation}
p'=\frac{4}{p}\textrm{, }\frac{p'}{2}-a'=\frac{\left(\frac{p}{2}-a\right)}{p}2i\sqrt{r}\textrm{, }q'=e^{-\pi/\sqrt{r}}.
\end{equation}
Also
\begin{equation}
Q_{\{a,p\}}(x)=\theta_{\{a,p\}}\left(e^{-\pi\sqrt{k_i(x)}}\right)
\end{equation}
and
\begin{equation}
y(x)=F\left(k\left(m_G(x)\right)\right)\textrm{, }k(r)=k_r.
\end{equation}
The function $F(x)$ is  $F(x)=R\left(e^{-\pi\sqrt{k_i(x)}}\right)$ and is a pure algebraic function i.e. it sends algebraic numbers to algebraic numbers and is evaluated as root of a sextic equation with rational coefficients (see [3],[15]). But $m_G(x)=k_i\left(Q_i(x)\right)$ and 
\begin{equation}
y(x)=F\left(Q_i(x)\right),
\end{equation}
where $Q(x)=Q_{\{a,p\}}(x)$, $Q_i(x)=Q^{(-1)}(x)$. From our conjectures $Q_{\{a,p\}}(x)$ is algebraic function. Hence $y(x)$ is also algebraic function. Inverting (54) we get
\begin{equation}
y_i(x)=\theta^{*}\circ k_i\circ F_i(x).
\end{equation}
From (49) we have the following\\ 
\\
\textbf{Theorem 5.}\\
If $\theta^{*}(r)=\theta_{\{a,p\}}(q)$ is as in (49), then 
\begin{equation}
y\left(\theta_{\{a,p\}}(q)\right)=R\left(q\right)\textrm{ and }m_G^{(-1)}(r)=\theta^{*}(r).
\end{equation}
\\

Also we have the next\\
\\
\textbf{Theorem 6.}\\
If $q=e^{-\pi\sqrt{r}}$, $r>0$, then\\
i) 
\begin{equation}
\frac{d\theta(q)}{dr}=\frac{1}{\phi(r)},
\end{equation}
where
\begin{equation}
\phi(r):=-\frac{2\sqrt{r}}{\pi \eta_D\left(\frac{i\sqrt{r}}{2}\right)^4G_{\{a,p\}}\left(R\left(q\right)\right)}
\end{equation}
and\\
ii)
\begin{equation}
\theta'(q)=q^{-5/6}\eta(q)^4G_{\{a,p\}}\left(R(q)\right).
\end{equation}
The function $G_{\{a,b\}}(x)$ is algebraic function of $x$.\\
\\
\textbf{Proof.}\\
From $m_G^{(-1)}(r)=\theta(q)$ we get 
\begin{equation}
5\int^{R(q)}_{0}\frac{G_{\{a,p\}}(t)dt}{t\sqrt[6]{t^{-5}-11-t^5}}=\theta(q).
\end{equation}
After derivating the above relation and using 
$$
\frac{dR(q)}{dq}=5^{-1}q^{-5/6}\eta(q)^4R(q)\sqrt[6]{R(q)^{-5}-11-R(q)^5},
$$
we get (ii):
$$
G_{\{a,p\}}\left(R(q)\right)q^{-5/6}\eta(q)^4=\theta'(q).
$$
Using the two definitions of Ramanujan's eta function and Dedekind's eta function, we get (i).\\
Of course $G(x)=G_{\{a,p\}}(x)$ is algebraic function and depends from $a,p$. From [15] we have
\begin{equation}
G_{\{a,p\}}\left(R(q)\right)=\sqrt[3]{\frac{k_r^2(1-k_r^2)}{2}}Q'_{\{a,p\}}(k_r).
\end{equation}
$qed$.\\
\\
We also have\\
\\
\textbf{Theorem 7.}\\
If $q=e^{-\pi\sqrt{r}}$, $r>0$, then
\begin{equation}
G_{\{a,p\}}(F(x))=\sqrt[3]{\frac{x^2(1-x^2)}{2}}Q'_{\{a,p\}}(x)\textrm{, }0<x<1
\end{equation}
and in the case of (49),(52) the function $G_{\{a,p\}}(x)$ is always algebraic.
\\

Now set in (59) where $q\rightarrow q^p$, then
$$
G\left(R\left(q^p\right)\right)q^{-5p/6}\eta\left(q^p\right)^4=\theta'\left(q^p\right).
$$
Hence
$$
\theta'\left(q^p\right)q^B=q^{B-5p/6}\eta\left(q^p\right)^4G\left(R\left(q^p\right)\right).
$$
Hence
$$
\sqrt[4]{\theta'\left(q^p\right)q^B}=q^{B/4-5p/24}\eta\left(q^p\right)\sqrt[4]{G\left(R\left(q^p\right)\right)}.
$$
Hence assuming that $B$ is such that $A=\frac{B}{4}-\frac{5p}{24}$, $A=-\frac{p}{12}+\frac{a}{2}-\frac{a^2}{2p}$, we have
$$
\sqrt[4]{\theta'\left(q^p\right)q^B}=\frac{\theta_4\left(\frac{p}{2},\frac{p}{2}-a;q\right)}{Q_{\{a,p\}}(k_r)}\sqrt[4]{G\left(R\left(q^p\right)\right)}.
$$
Hence using (52) we get:\\
\\
\textbf{Theorem 8.}\\If $|q|<1$ and $B=4A+5p/6$, $A=-\frac{p}{12}+\frac{a}{2}-\frac{a^2}{2p}$, then
$$
\theta'\left(q^p\right)=q^{-B}\left(\frac{G\left(R\left(q^p\right)\right)}{\theta\left(q\right)^4}\right)\theta_4\left(\frac{p}{2},\frac{p}{2}-a;q\right)^4=
$$
\begin{equation}
=q^{-B}\sqrt[3]{\frac{k_{p^2r}^2(1-k_{p^2r}^2)}{2}}\frac{Q_{\{a,p\}}'(k_{p^2r})}{(Q_{\{a,p\}}(k_r))^4}\theta_4\left(\frac{p}{2},\frac{p}{2}-a;q\right)^4.
\end{equation}
\\
\textbf{Example 1.}\\
Suppose $a=1$, $p=4$, then
\begin{equation}
\sum^{\infty}_{n=-\infty}(-1)^nq^{2n^2+n}=q^{1/24} \eta\left(q^4\right)Q_{\{1,4\}}(k_r)
\end{equation}
Then $Q_{\{1,4\}}(x)$ will be 
$$
Q_{\{1,4\}}(x)=\sqrt[12]{4\frac{1-x^2}{x}}
$$
For a certain $G$ we have (for the function $\sigma(x)$ see [15]): 
$$
\sigma(x)=-6\cdot2^{5/6}x^{13/12}(1-x^2)^{11/12}(1+x^2)^{-1}.
$$  
Hence
$$
G_{\{1,4\}}\left(R(q)\right)=-\frac{1+k_r^2}{12\cdot 2^{1/6}k_r^{5/12}(k^{*}_r)^{7/6}}\textrm{, }k^{*}_r=\sqrt{1-k_r^2}.
$$
and
$$
y(x)=R\left(e^{-\pi\sqrt{k_i\left(\frac{1}{8}\left(-x^{12}+\sqrt{64+x^{24}}\right)\right)}}\right).
$$
\\
\textbf{Example 2.}\\
From relation (14), whenever $s$ is integer we have 
\begin{equation}
\theta_3(1,2s+1;q)=\sum^{\infty}_{n=-\infty}q^{n^2+(2s+1)n}
=q^{-1/6-s-s^2}\eta(q^2)\frac{2^{1/3}\sqrt[3]{k}}{\sqrt[6]{k^{*}}}
\end{equation}
and
\begin{equation}
Q_{\{-2s,2\}}(x)=\frac{2^{1/3}x^{1/3}}{(1-x^2)^{1/12}}.
\end{equation}
Hence
\begin{equation}
\theta_{\{-2s,2\}}'\left(q^2\right)=\frac{q^{4s^2+4s-1}}{24\cdot 2^{5/6}}\frac{1+6k^{*}+(k^{*})^2}{(1-k^*)^{2/3}(k^{*})^{1/12}(1+k^*)^{7/6}}\cdot \theta^4_3(1,2s+1;q),
\end{equation}
\begin{equation}
G\left(R(q)\right)=\frac{2-k_r^2}{6(k^{*}_r)^{3/2}}
\end{equation}
and
\begin{equation}
\theta'(q)=\theta_{\{-2s,2\}}'(q)=q^{-5/6}\eta(q)^4\frac{2-k_r^2}{6(k^{*}_r)^{3/2}}
\end{equation}
and
\begin{equation}
y(x)=R\left(\exp\left[-\pi\sqrt{k_i\left(\frac{1}{4\sqrt{2}}\sqrt{-x^{12}+x^6\sqrt{64+x^{12}}}\right)}\right]\right).
\end{equation}
\\
\textbf{Theorem 9.}\\
Assume the DE 
\begin{equation}
X'(x)+\frac{2^{4/3}}{\left(x\sqrt{1-x^2}\right)^{2/3}}P\left(X(x)\right)=0.
\end{equation}
Then if
\begin{equation}
Y(r)=X(k_r)
\end{equation}
we have
\begin{equation}
Y'(r)=\frac{dY(r)}{dr}=\pi\frac{\eta_D\left(i\sqrt{r}/2\right)^4}{\sqrt{r}}P\left(Y(r)\right).
\end{equation}
\\

Hence if we assume that exists $P(x)$ as above and
\begin{equation}
\theta_{\{a,p\}}(q)=Q_{\{a,p\}}(k_r)=Y(r)=X(k_r),
\end{equation}
then
\begin{equation}
Q_{\{a,p\}}(k_r)=\theta_{\{a,p\}}(q)\textrm{, }q=e^{-\pi\sqrt{r}}\textrm{, }r>0.
\end{equation}
But from Theorems 7,9 we have
$$
G_{\{a,p\}}\left(R(q)\right)=\frac{\left(k_rk^{*}_r\right)^{3/2}}{2^{1/3}}Q'_{\{a,p\}}(k_r)=
\frac{\left(k_rk^{*}_r\right)^{2/3}}{2^{1/3}}\left(\frac{-2^{4/3}P\left(X(k_r)\right)}{\left(k_rk^{*}_r\right)^{2/3}}\right)=
$$
$$
=-2P\left(Q_{\{a,p\}}(k_r)\right).
$$
Hence we get the next\\
\\
\textbf{Theorem 10.}\\
If $r>0$, then
\begin{equation}
G_{\{a,p\}}\left(R(q)\right)=-2P\left(Q_{\{a,p\}}(k_r)\right)
\end{equation}
and
\begin{equation}
P\left(Q_{\{a,p\}}(x)\right)=-2^{-4/3}\frac{\left(x\sqrt{1-x^2}\right)^{2/3}}{\sigma(x)}=-\frac{x^{2/3}(\sqrt{1-x^2})^{2/3}}{2^{4/3}}Q'_{\{a,p\}}(x),
\end{equation}
which is esentialy (71).\\
\\
\textbf{Lemma.}\\
If $r>0$, then
\begin{equation}
-\frac{\pi}{2}\frac{\eta_D\left(i\sqrt{r}/2\right)^4}{\sqrt{r}}=\frac{d}{dr}\left(\frac{1}{\sqrt[3]{4}}B\left(k_r^2;\frac{1}{6},\frac{2}{3}\right)\right).
\end{equation}
\\
\textbf{Proof.}\\
See Theorem 17 below.\\
\\

Using the above Lemma and integrating (73) we get
\begin{equation}
-\int^{Q_{\{a,p\}}(y)}_{Q_{\{a,p\}}(x)}\frac{dt}{P(t)}=\frac{2}{\sqrt[3]{4}}B\left(y^2;\frac{1}{6},\frac{2}{3}\right)-\frac{2}{\sqrt[3]{4}}B\left(x^2;\frac{1}{6},\frac{2}{3}\right),
\end{equation}
where $x,y\in (0,1)$.\\
\\
\textbf{Theorem 11.}\\
For the function $Y(r)$ of (72),(73) it holds
\begin{equation}
-\frac{1}{2}\int^{Y(r)}_{Y(\infty)}\frac{dt}{P(t)}=\frac{1}{\sqrt[3]{4}}B\left(k_r^2;\frac{1}{6},\frac{2}{3}\right)\textrm{, }r>0.
\end{equation}
Also in the special case of a theta function, we have
$$
Y(r)=Q_{\{a,p\}}(k_r).\eqno{(80.1)}
$$
\\
\textbf{Proof.}\\
Easy from Lemma and (73).\\
\\
\textbf{Theorem 12.}\\
If $A$ is real and $X(0)=0$ and 
\begin{equation}
X(A)=\int^{A}_{0}\frac{dt}{\sigma(t)}=s_i(A)\left(=Q_{\{a,p\}}(A)\right),
\end{equation}
then $P(A)$ is such that
\begin{equation}
P(A)=-\frac{1}{2h_i'(A)}\Leftrightarrow -\frac{1}{2}\int^{h(A)}_{c}\frac{dt}{P(t)}=A.
\end{equation}
Hence given $P(A)$, the solution of (71) is
\begin{equation}
X(A)=h\left(\frac{1}{\sqrt[3]{4}}B\left(A^2;\frac{1}{6},\frac{2}{3}\right)\right),
\end{equation}
where $h(A)$ is that of (82).
Also given a function $X(A)$ such  
\begin{equation}
X(A)=h\left(\frac{1}{\sqrt[3]{4}}B\left(A^2;\frac{1}{6},\frac{2}{3}\right)\right),
\end{equation} 
then $X(A)$ is that (81) and $P(A)$ is that of (82) and (84) satisfies (71).\\
\\
\textbf{Proof.}\\
Assume $X(x)=s_i(x)$, then from Theorem 9 we have
$$
s_i'(A)+\frac{2^{4/3}}{\left(A\sqrt{1-A^2}\right)^{2/3}}P\left(s_i(A)\right)=0,
$$
thus
$$
\frac{1}{s'(A)}+\frac{2^{4/3}}{\left(s(A)\sqrt{1-s(A)^2}\right)^{2/3}}P(A)=0,
$$
thus
$$
\frac{1}{P(A)}=-\frac{2^{4/3}s'(A)}{\left(s(A)\sqrt{1-s(A)^2}\right)^{2/3}}.
$$
But
\begin{equation}
h_i'(A)=\frac{2^{1/3}s'(A)}{\left(s(A)\sqrt{1-s(A)^2}\right)^{2/3}}.
\end{equation}
Hence
$$
P(A)=\frac{1}{-2h_i'(A)}.
$$
Also from
\begin{equation}
h\left(\frac{1}{\sqrt[3]{4}}B\left(A^2;\frac{1}{6},\frac{2}{3}\right)\right)=s_i(A),
\end{equation}
we get (82) and (81).\\
\\
\textbf{Notes.}\\
Theorems 9,11,12 are in accordance with each other. Also in case $X(A)=Q_{\{a,p\}}(A)$, then
$$
h(A)=Q_{\{a,p\}}\left(k\left(m(A)\right)\right),\eqno{(86.1)}
$$
where $k(A)=k_A$. The function $m(A)$ is defined as the function which satisfies 
$$
\pi\int^{+\infty}_{\sqrt{m(r)}}\eta_D(it/2)^4dt=r\eqno{(86.2)}
$$
and 
$$
m\left(\frac{1}{\sqrt[3]{4}}B\left(k_r^2;\frac{1}{6},\frac{2}{3}\right)\right)=r.\eqno{(86.3)}
$$
Hence we have the next:\\
\\ 
\textbf{Theorem 13.}\\
For every function $Y=Y(r)$, we define $X(x)$ as $Y(r)=X(k_r)$. Then  $X(x)$ is solution of a DE 
$$
X'(x)+\frac{2^{4/3}}{\left(x\sqrt{1-x^2}\right)^{2/3}}P\left(X(x)\right)=0.
$$
If also $\lim_{t\rightarrow 0}X(t)=c$, then the equation $Y(r)=Y_0$ have solution
$$
r=m\left(-\frac{1}{2}\int^{Y_0}_{c}\frac{dt}{P(t)}\right).\eqno{(86.4)}
$$
Also holds
$$
F_c\left(Y\left(4 r\right)\right)+F_c\left(Y\left(\frac{4}{r}\right)\right)=-\frac{\sqrt{3}\Gamma\left(\frac{1}{3}\right)^3}{\pi\sqrt[3]{2}}\textrm{, }\forall r>0,\eqno{(86.4.1)}
$$
where
 $$
 F_c(x)=\int^{x}_{c}\frac{dt}{P(t)}.
 $$
\\
\textbf{Proof.}\\
If $X(x)$ is a solution of DE (71) and $Y(r)=X(k_r)$, $c=\lim_{t\rightarrow0}X(t)$, then from Theorem 11 and definition (86.2),(86.3), we have the result.\\
\\
\textbf{Notes.}\\
\textbf{i)} If we set
$$
F_c(x):=\int^{x}_{c}\frac{dt}{P(t)},\eqno{(86.5)}
$$
and
$$
B(x):=3\sqrt[3]{2x}\cdot{}_2F_1\left(\frac{1}{3},\frac{1}{6};\frac{7}{6};x^2\right)=\frac{1}{\sqrt[3]{4}}B\left(x^2;\frac{1}{6},\frac{2}{3}\right),\eqno{(86.6)}
$$
then
\begin{equation}
B\left(k_{4r}\right)+B\left(k_{4/r}\right)=\frac{\sqrt{3}\Gamma\left(\frac{1}{3}\right)^3}{2\pi \sqrt[3]{2}},
\end{equation}
and we get
\begin{equation}
F_c\left(Y\left(4 r\right)\right)+F_c\left(Y\left(\frac{4}{r}\right)\right)=-\frac{\sqrt{3}\Gamma\left(\frac{1}{3}\right)^3}{\pi\sqrt[3]{2}}\textrm{, }\forall r>0.
\end{equation}
In the special case which $\theta_{\{a,p\}}(q)=Q_{\{a,p\}}(k_r)=Q(k_r)$, then $Y(r)=\theta^{*}(r)=Q_{\{a,p\}}(k_r)$.\\
\textbf{ii)} Assuming that every algebraic function $U(x)$ is again algebraic function of the singular modulus $k_{r}$ i.e. $U(x)=X(k_r)=Y(r)$ (this hapens, since for every $x$ algebraic, exists algebraic functions $X^{(-1)}(t)$ and $X_1(t)$ such that $X_1(k_r)=x$ and $X^{(-1)}(U(X_1(t)))=t$). Then exists a new algebraic function $P(x)$ such that   
\begin{equation}
\frac{dU(x)}{dk}=\frac{2^{4/3}}{(kk')^{2/3}}P\left(U(x)\right).
\end{equation}
Hence 
$$
\frac{dU(x)}{dr}=X'(k_r)\frac{dk_r}{dr}=Y'(r)=\pi\frac{\eta_D\left(i\sqrt{r}/2\right)^4}{\sqrt{r}}P\left(Y(r)\right),
$$
from
$$
\frac{dk_r}{dr}=\frac{k_r(k'_r)^2K^2}{\pi\sqrt{r}}
$$
and (9.1) and (89). Hence
$$
\frac{dU(x)}{dr}=\frac{2^{4/3}k^{1/3}(k')^{4/3}K^2}{\pi\sqrt{r}}P\left(U(x)\right).
$$
By this way for every algebraic function $U(x)=Y(r)$, there exists algebraic function $P(x)$ such that
\begin{equation}   
Y'(r)=\pi\frac{\eta_D\left(i\sqrt{r}/2\right)^4}{\sqrt{r}}P\left(Y(r)\right).  
\end{equation}
The function $Y(r)$ is Hauptmodul according to (88). Hence\\  
\\
\textbf{Theorem 13.1.}\\ 
Every algebraic function $U(x)$ can be writen in the form $U(x)=X(k_r)=Y(r)$. The function $Y(r)$ is a Hauptmodul and satisfies (88) and (90).\\
\\  

We return again to Theorem 12 and we combine Theorem 6 of [15] with the function $P(x)$. It is already known, that $P$ and $h$ are related with (82) and if 
\begin{equation}
h_1(t):=\left(\frac{1}{h_i'(x)}\right)^{(-1)}(t),
\end{equation}
then we have\\
\\
\textbf{Theorem 14.}(see [15]) 
\begin{equation}
5\int^{G_i(x)}_{0}\frac{dt}{t\sqrt[6]{t^{-5}-11-t^5}}=\int^{x}_{c}\frac{h'_1(t)}{t}dt.
\end{equation}
Moreover
\begin{equation}
G_i(x)=F_1\left(\int^{x}_{c}\frac{h_1'(t)}{t}dt\right),
\end{equation}
where $F_1(x)$ is defined from
\begin{equation}
x=5\int^{F_1(x)}_{0}\frac{dt}{t\sqrt[6]{t^{-5}-11-t^5}}.
\end{equation}
\\

Hence easily from the above 
\begin{equation}
h_1\left(-2P(A)\right)=A\Rightarrow -2h_1'(-2P(A))P'(A)=1
\end{equation}
and 
$$
\int^{x}_{c}\frac{h_1'(t)}{t}dt=\int^{P_i\left(-x/2\right)}_{c_1}\frac{h_1'(-2P(A))}{-2P(A)}(-2P'(A))dA=
$$
$$
=\int^{P_i(-x/2)}_{c_1}\frac{1}{4P'(A)P(A)}(-2P'(A))dA=-\frac{1}{2}\int^{P_i(-x/2)}_{c_1}\frac{dA}{P(A)}.
$$
Hence we can state the next\\ 
\\
\textbf{Theorem 15.}
\begin{equation}
G^{(-1)}(-2x)=F_1\left(\frac{-1}{2}\int^{P_i(x)}_{c}\frac{dt}{P(t)}\right).
\end{equation}
\\
\textbf{Application.}\\
If $q=e^{-\pi\sqrt{r}}$, $r>0$, then  $\theta_3(q)=\sum^{\infty}_{n=-\infty}q^{n^2}=\theta_3(1,0;q)$. Hence from (9.1) and $\theta_3(q)=\sqrt{2K(k_r)/\pi}$, we have
$$
\theta^{*}(r)=\theta_{\{1,2\}}(q)=A_{\{1,2\}}(k_r)=Q_{\{1,2\}}(k_r)=q^{-1/12}\frac{\theta_3(q)}{\eta(q^2)}=
$$
$$
=\frac{2^{1/6}}{(k_{4r})^{1/2}(k^{*}_{4r})^{1/3}}\sqrt{\frac{K(k_r)}{K(k_{4r})}}.
$$
But
$$
k_{4r}=\frac{1-k^{*}_r}{1+k^{*}_r}\textrm{ and }\frac{K(k_{4r})}{K(k_r)}=\frac{1+k^{*}_r}{2},
$$
where $k^{*}=\sqrt{1-k^2}$. Hence
$$
Q_{\{1,2\}}(k)=\frac{\sqrt[3]{2}}{k^{1/6}(1-k^2)^{1/12}}.
$$
From Theorem 12 we have
$$
X(A)=\frac{\sqrt[3]{2}}{A^{1/6}(1-A^2)^{1/12}},
$$
then 
$$
P\left(\frac{\sqrt[3]{2}}{x^{1/6}(1-x^2)^{1/12}}\right)=\frac{1-2x^2}{12\sqrt{x}(1-x^2)^{3/4}}.
$$
Also
$$
\sigma\left(x\right)=\frac{3\cdot 2^{2/3}x^{7/6}(1-x^2)^{13/12}}{-1+2x^2}.
$$
\\

By this way one can see that, after we define equation (71) of Theorem 9, all functions of [15] are become meaningfull and evaluated easily. Next we shall examine equation (71) and answer the question: from where (71) comes and what it represents? For to answer this, we shall give generalizations of some functions used so far in the complex plane.

\section{The Complex Analog of Theorem 9 and Hauptmodul's}

We re-define
\begin{equation}
K(w):=\frac{\pi}{2}{}_2F_1\left(\frac{1}{2},\frac{1}{2};1;w^2\right)\textrm{, }|w|<1,
\end{equation}
to be the complete elliptic integral of the first kind, assuming that takes and complex values. Assume also $q=e_1(z):=e^{\pi i z}$, with $Im(z)>0$ and  
\begin{equation}
m^{*}(z):=\left(\frac{\theta_2\left(e^{i\pi z}\right)}{\theta_3\left(e^{i\pi z}\right)}\right)^2.
\end{equation}
Then it is known that if $z=x+iy$, $-\frac{1}{2}< x\leq  \frac{1}{2}$, $y>0$, we have
\begin{equation}
i\frac{K\left(\sqrt{1-m^{*}(z)^2}\right)}{K\left(m^{*}(z)\right)}=z.
\end{equation}
We also set the complete elliptic integral of the first kind at singular values to be 
\begin{equation}
K:=K[z]:=K(m^{*}(z))=\frac{\pi}{2}\cdot {}_2F_1\left(\frac{1}{2},\frac{1}{2};1;m^{*}(z)^2\right).
\end{equation}
The Dedekind eta function $\eta_D(z)$, can evaluated by means of the singular modulus $m^{*}(z)$ and $K(m^{*}(z))$ using the next formula:
$$
\eta_D(z)^4=\frac{2^{4/3}}{\pi^2}m^{*}(2z)^{1/3}\left(1-m^{*}(2z)^2\right)^{2/3}K[2z]^2=
$$
\begin{equation}
=\frac{2^{2/3}}{\pi^2}m^{*}(z)^{2/3}\left(1-m^{*}(z)^2\right)^{1/3}K[z]^2,
\end{equation}
which is similar to (9.1), but now is defined in the complex upper half plane.\\
It is also well known that $\eta_D(z)$ have modular properties, since
\begin{equation}
\eta_D\left(-\frac{1}{z}\right)=\sqrt{-iz}\cdot\eta_D(z)\textrm{, }Im(z)>0.
\end{equation}
We also give the definition of Hauptmodul functions: A  function $f(z)$, defined in the upper half plane $H$ is called Hauptmodul, if exist a function $g$ such that
\begin{equation}
f\left(-\frac{1}{z}\right)=g(f(z))\textrm{, }\forall z\in H.
\end{equation}
Such functions are the Klein's $j-$invariant i.e.
\begin{equation}
j(z)=\frac{E_4(z)^3}{\Delta(z)}=q^{-1}+744+196884q+\ldots,
\end{equation}
where $E_4(z)$ is the Eisenstein weight 4 modular form ($q=e(z)$):
\begin{equation}
E_4(z)=1+240\sum^{\infty}_{n=1}\sigma_{3}(n)q^n\textrm{, }\sigma_3(n)=\sum_{d|n}d^3
\end{equation}
and $\Delta(z)$ (a cusp form)
\begin{equation}
\Delta(z)=q\prod^{\infty}_{n=1}(1-q^n)^{24}=\eta_D(z)^{24}.
\end{equation}
Then
\begin{equation}
j\left(-\frac{1}{z}\right)=j(z)\textrm{, }Im(z)>0.
\end{equation}
Also if $R^{*}(z)$ denotes the Rogers-Ramanujan continued fraction for the argument $z$ in $q=e(z):=e^{2\pi i z}$:
\begin{equation}
R^{*}(z)=\frac{q^{1/5}}{1+}\frac{q^1}{1+}\frac{q^2}{1+}\frac{q^3}{1+}\ldots=q^{1/5}\prod^{\infty}_{n=1}(1-q^n)^{(n|5)},
\end{equation}
where $(n|l)$ is the well known Jacobi symbol, then we have
\begin{equation}
R^{*}\left(-\frac{1}{z}\right)=\frac{1-\phi R^{*}(z)}{\phi+R^{*}(z)}\textrm{, }Im(z)>0
\end{equation}
and $\phi=\frac{1+\sqrt{5}}{2}$ is the golden ratio.\\
Another example is Carty's function $\Pi(z)$ (see Theorems 17,18,19 below). The function $\Pi(z)$ is defined when $-\frac{1}{2}<Re(z)\leq \frac{1}{2}$ and $Im(z)>0$ as
\begin{equation}
\Pi(z):=-2\pi i\int^{i\infty}_{z}\eta(t)^4dt=\frac{1}{\sqrt[3]{4}}B\left(m^{*}(2z)^2;\frac{1}{6},\frac{2}{3}\right).
\end{equation}
Then also
\begin{equation}
\Pi\left(-\frac{1}{z}\right)=\frac{\sqrt{3}\Gamma\left(\frac{1}{3}\right)^3}{2\pi \sqrt[3]{2}}-\Pi(z)\textrm{, }Im(z)>0.
\end{equation}
Our final example and most common is the singular modulus itself
\begin{equation}
m^{*}\left(-\frac{1}{z}\right)=\sqrt{1-m^{*}(z)^2},
\end{equation}
where also $-\frac{1}{2}<Re(z)\leq \frac{1}{2}$ and $Im(z)>0$.\\
One can naturaly ask: are there more such examples, and further, is there a convinient and easy way to construct such functions? The answer we give here is afirmative.\\
\\	

Assume the differential equations 
\begin{equation}
Y'(z)+4\pi i\cdot \eta_D(z)^4P\left(Y(z)\right)=0
\end{equation}
and
\begin{equation}
Y'(z)-4\pi i\cdot \eta_D(z)^4P(Y(z))=0.
\end{equation}
These can be unified if we write them as
\begin{equation}
Y'(z)^2+16\pi^2\cdot \eta_D(z)^8 P\left(Y(z)\right)^2=0.
\end{equation}
We have a first result.\\
\\
\textbf{Theorem 16.}\\
Given a smooth function $P(x)$, $x\in\textbf{R}$, we consider the differential equation
\begin{equation}
X'(x)+2^{4/3}x^{-2/3}(1-x^2)^{-1/3}P\left(X(x)\right)=0
\end{equation}
and set $Y(z)=X(m^{*}(2z))$, where $X(x)$ is solution of (116). Then $Y_1(z)=Y(z)$, $Y_2(z)=Y(-1/z)$, with $Im(z)>0$ are solutions of (115). More precicely $Y_1(z)$ is solution of (113) and $Y_2(z)$ is solution of (114).\\
\\
\textbf{Proof.}\\
If we set $Y(z)=X(m^{*}(2z))$, where $X(x)$ is solution of
$$
X'(x)+2^{4/3}x^{-2/3}(1-x^2)^{-1/3}P\left(X(x)\right)=0,\eqno{(eq)}
$$ 
then we have from [9],[20]:
\begin{equation}
{m^{*}}'(z)=\frac{dm^{*}(z)}{dz}=\frac{2i}{\pi}m^{*}(z)(1-m^{*}(z)^2)K^2.
\end{equation}
Hence
$$
Y_1'(z)=2X'\left(m^{*}(2z)\right)\frac{dm^{*}(2z)}{dz}=
$$
\begin{equation}
=\frac{4i}{\pi}m^{*}(2z)\left(1-m^{*}(2z)^2\right)K[2z]^2X'\left(m^{*}(2z)\right).
\end{equation}
Substituting in $(eq)$ where $x\rightarrow m^{*}(2z)$, we get
\begin{equation}
X'(m^{*}(2z))+2^{4/3}m^{*}(2z)^{-2/3}\left(1-m^{*}(2z)^2\right)^{-1/3}P\left(X(m^{*}(2z))\right)=0.
\end{equation}
If we multiply both sides of the above equation with ${m^{*}}'(2z)$, we get     
$$
X'\left(m^{*}(2z)\right){m^{*}}'(2z)=
$$
$$
=-2^{4/3}{m^{*}}'(2z)m^{*}(2z)^{-2/3}\left(1-m^{*}(2z)^2\right)^{-1/3}P\left(X(m^{*}(2z))\right).
$$
Hence using (117), we get
$$
Y_1'(z)+2^{4/3}\frac{2i}{\pi}\frac{m^{*}(2z)\left(1-m^{*}(2z)^2\right)}{m^{*}(2z)^{2/3}\left(1-m^{*}(2z)^2\right)^{1/3}}K[2z]^2P\left(Y_1(z)\right)\Leftrightarrow
$$
$$
Y_1'(z)+2^{4/3}\frac{2i}{\pi}m^{*}(2z)^{1/3}\left(1-m^{*}(2z)^2\right)^{2/3}K[2z]^2P\left(Y_1(z)\right)=0.
$$
Using (101) we get the result for the first equation
$$
Y_1'(z)+4\pi i \eta_D(z)^4P\left(Y_1(z)\right)=0.
$$
Now we prove the result for the $Y_2(z)$ function.\\
We have 
\begin{equation}
Y_2'(z)=X'\left(m^{*}\left(-\frac{2}{z}\right)\right){m^{*}}'\left(-\frac{2}{z}\right)\frac{1}{z^2}.
\end{equation}
Also if we set $x\rightarrow m^{*}\left(-\frac{2}{z}\right)$ in $(eq)$, we get
$$
X'\left(m^{*}\left(\frac{-2}{z}\right)\right)+
$$
\begin{equation}
+2^{4/3}m^{*}\left(\frac{-2}{z}\right)^{-2/3}\left(1-m^{*}\left(\frac{-2}{z}\right)^2\right)^{-1/3}P\left(X\left(m^{*}\left(\frac{-2}{z}\right)\right)\right)=0
\end{equation}
But it is also known that 
\begin{equation}
m^{*}\left(\frac{-1}{z}\right)^2=1-m^{*}(z)^2
\end{equation}
and 
\begin{equation}
m^{*}\left(2z\right)=\frac{1-\sqrt{1-m^{*}(z)^2}}{1+\sqrt{1-m^{*}(z)^2}}.
\end{equation}   
Hence
\begin{equation}
m^{*}\left(\frac{-2}{z}\right)=\frac{1-m^{*}(z)}{1+m^{*}(z)}.
\end{equation}
Also
\begin{equation}
K\left(m^{*}\left(\frac{-2}{z}\right)\right)=-iK\left(m^{*}\left(\frac{z}{2}\right)\right)\frac{z}{2}.
\end{equation}
From (117) we can write  
\begin{equation}
{m^{*}}'\left(\frac{-2}{z}\right)=\frac{2i}{\pi}m^{*}\left(\frac{-2}{z}\right)\left(1-m^{*}\left(\frac{-2}{z}\right)^2\right)K\left(m^{*}\left(\frac{-2}{z}\right)\right)^2.
\end{equation} 
Hence we get
$$
-X'\left(m^{*}\left(-\frac{2}{z}\right)\right){m^{*}}'\left(-\frac{2}{z}\right)\frac{1}{z^2}=
$$
$$
2^{4/3}{m^{*}}'\left(\frac{-2}{z}\right)\frac{2}{z^2}m^{*}\left(\frac{-2}{z}\right)^{-2/3}\left(1-m^{*}\left(\frac{-2}{z}\right)^2\right)^{-1/3}P\left(X\left(m^{*}\left(\frac{-2}{z}\right)\right)\right).
$$ 
Using (124),(125),(126), we get after simplifications
$$
Y_2'(z)-\frac{2i2^{2/3}\left(1-m^{*}(z)\right)^{1/3}m^{*}(z)^{2/3}}{\pi(1+m^{*}(z))^{5/3}}K\left[z/2\right]^2P\left(Y_2(z)\right)=0.
$$
Using now
\begin{equation}
K[z/2]=(1+m^{*}(z))K[z],
\end{equation}
we get
$$
Y_2'(z)=\frac{2i2^{2/3}\left(1-m^{*}(z)\right)^{1/3}m^{*}(z)^{2/3}}{\pi(1+m^{*}(z))^{5/3}}(1+m^{*}(z))^2K\left[z\right]^2P\left(Y_2(z)\right)\Leftrightarrow
$$
$$
Y_2'(z)=4\pi i\frac{2^{2/3}}{\pi^2}\left(1-m^{*}(z)^2\right)^{1/3}m^{*}(z)^{2/3}K\left[z\right]^2P(Y_2(z))\Leftrightarrow
$$
$$
Y_2'(z)=4\pi i \eta_D(z)^4P\left(Y_2(z)\right)
$$
which is the second equation of (115). These complete the  validity of the theorem.\\
\\  

Next we state a very important theorem about the generalized integral of $\eta_D(z)^4$.\\
\\
\textbf{Theorem 17.}(Carty's theorem)\\
Whenever $Im(w)>0$ and $Im(z)>0$, we have
\begin{equation}
2\pi i\int^{w}_{z}\eta_D(t)^4dt=\left[\frac{1}{\sqrt[3]{4}}B\left(m^{*}(2t)^2;\frac{1}{6},\frac{2}{3}\right)\right]^{t=w}_{t=z},
\end{equation}
where $B(z;a,b):=\int^{z}_{0}t^{a-1}(1-t)^{b-1}dt$ is the incomplete Beta function.\\ 
\\
\textbf{Proof.}\\
For to prove (128), we derivate with respect to $w$, both sides of the idenity and get easily
$$
2\pi i \eta_D(w)^4=2^{4/3}\frac{{m^{*}}'(2w)}{\left(1-m^{*}(2w)^2\right)^{1/3}m^{*}(2w)^{2/3}}.
$$ 
Now we use (117) to arive at
$$
2\pi i\eta_D(w)^4=
$$
$$
=\frac{2^{4/3}}{(1-m^{*}(2w)^2)^{1/3}m^{*}(2w)^{2/3}}\frac{2i}{\pi}m^{*}(2w)\left(1-m^{*}(2w)^2\right)K\left(m^{*}(2w)\right)^2.
$$
Then after simplification we arrive to
$$
\eta_D(w)^4=\frac{2^{4/3}}{\pi^2}m^{*}(2w)^{1/3}\left(1-m^{*}(2w)^2\right)^{2/3}K[2w]^2,
$$
which is true (see relation (101)). Hence going backwords we get easily the result.\\
\\

At this point we can say that we have one equation
$$
(Y'(z))^2+16\pi^2\eta_D(z)^8P^2\left(Y(z)\right)=0,
$$ 
which have solutions $Y_1(z)=X(m^{*}(2z))$, $Y_2(z)=X\left(m^{*}(-2/z)\right)$ and can be writen as
$$
Y'(z)=\pm 4\pi i \eta_D(z)^4P(Y(z))
$$
Hence we can integrate this last equation using Carty's theorem to get
\begin{equation}
\int^{Y_1(w)}_{Y_1(z)}\frac{dt}{P(t)}=-\left[\frac{2}{\sqrt[3]{4}}B\left(m^{*}(2t)^2;\frac{1}{6},\frac{2}{3}\right)\right]^{t=w}_{t=z}
\end{equation}
and
\begin{equation}
\int^{Y_1(-1/w)}_{Y_1(-1/z)}\frac{dt}{P(t)}=\left[\frac{2}{\sqrt[3]{4}}B\left(m^{*}(2t)^2;\frac{1}{6},\frac{2}{3}\right)\right]^{t=w}_{t=z}.
\end{equation}
Hence adding both relations:
\begin{equation}
\int^{Y_1(w)}_{Y_1(z)}\frac{dt}{P(t)}+\int^{Y_1(-1/w)}_{Y_1(-1/z)}\frac{dt}{P(t)}=0
\end{equation}
or equivalently
$$
\int^{Y_1(w)}_{c}\frac{dt}{P(t)}+\int^{Y_1(-1/w)}_{c}\frac{dt}{P(t)}=\int^{Y_1(z)}_{c}\frac{dt}{P(t)}+\int^{Y_1(-1/z)}_{c}\frac{dt}{P(t)}=C_1,
$$
where $c$ and $C_1$ are constants. Set now $c=Y_1(z_0)$ and $c'=Y_1(-1/z_0)$, with $Im(z_0)>0$. When $Im(z)>0$, we have
\begin{equation}
\int^{Y_1(z)}_{c'}\frac{dt}{P(t)}+\int^{Y_1(-1/z)}_{c}\frac{dt}{P(t)}=C_1.
\end{equation}
Define now the function $F_c(x)$ as
\begin{equation}
F_c(x):=F(c,x):=\int^{x}_{c}\frac{dt}{P(t)},
\end{equation}
with $c$ constant. From (131),(132) we finaly get
\begin{equation}
Y_1\left(-\frac{1}{z}\right)=F^{(-1)}_{c'}\left(C_1\pm F_c\left(Y_1(z)\right)\right)\textrm{, }\forall z\in H_1\subset H
\end{equation}
From the above notes it is clear that we can state the next\\
\\ 
\textbf{Theorem 18.}\\
Let $P(x)$ be a ''suitable'' smooth function in a local interval and $X(x)$ a function defined from the following DE
\begin{equation}
X'(x)+2^{4/3}x^{-2/3}(1-x^2)^{-1/3}P\left(X(x)\right)=0.
\end{equation}
If 
$$
Y(z)=X\left(m^{*}(2z)\right)\textrm{, }z\in H_1\subset H
$$
and 
$$
F_c(x)=F(c,x)=\int^{x}_{c}\frac{dt}{P(t)},
$$
with $c$ constant, then
$$
Y\left(-\frac{1}{z}\right)=F^{(-1)}_{c'}\left(C_1\pm F_c\left(Y(z)\right)\right)\textrm{, }\forall z\in H_1\subset H.\eqno(135.1)
$$
The constants $c$, $c'$ are such that for any fixed $z_0\in H$, $c=Y(z_0)$ and $c'=Y(-1/z_0)$.\\
\\
\textbf{Notes.}\\
Theorem 18 says that for every smooth function $Y(z)$, under some weak conditions (here $P(t)\neq 0$ in a certain region), there exists a function $F_c(z)$ such that (135.1) holds for all $z\in H_1$. Hence in a local region every smooth function is a Hauptmodul.\\
\\   
\textbf{Theorem 19.}\\
The function $Y(z)=X\left(m^{*}(2z)\right)$ is defined from the equation
\begin{equation}
\int^{Y(z)}_{Y(w)}\frac{dt}{P(t)}=-\frac{2}{\sqrt[3]{4}}\left(B\left(m^{*}(2z)^2;\frac{1}{6},\frac{2}{3}\right)-B\left(m^{*}(2w)^2;\frac{1}{6},\frac{2}{3}\right)\right)
\end{equation}
and the oposite.
Also if $c=Y(i\infty)$ and $F_c(x)$ as in Theorem 18, we get 
$$
\int^{Y(z)}_{c}\frac{dt}{P(t)}=-\frac{2}{\sqrt[3]{2}}B\left(m^*(2z)^2;\frac{1}{6},\frac{2}{3}\right)\textrm{, }Im(z)>0\eqno{(136.1)}
$$
\begin{equation}
Y\left(-\frac{1}{z}\right)=F^{(-1)}_c\left(\frac{-\sqrt{3}\Gamma\left(\frac{1}{3}\right)^3}{\pi \sqrt[3]{2}}-F_c\left(Y(z)\right)\right)\textrm{, }H_1\subset H.
\end{equation}
\\
\textbf{Proof.}\\
Set 
$$
\Pi_1(z):=\frac{2}{\sqrt[3]{4}}B\left(m^{*}(2z);\frac{1}{6},\frac{2}{3}\right)\textrm{, }z\in H.
$$
It holds from Theorems 16,17 and relation (113):
\begin{equation}
\int^{Y(z)}_{Y(i\infty)}\frac{dt}{P(t)}=-\frac{2}{\sqrt[3]{4}}B\left(m^{*}(2z)^2;\frac{1}{6},\frac{2}{3}\right)
\end{equation}
and
\begin{equation}
\int^{Y(-1/z)}_{Y(i\infty)}\frac{dt}{P(t)}=-\frac{2}{\sqrt[3]{4}}B\left(m^{*}\left(-\frac{2}{z}\right)^2,\frac{1}{6};\frac{2}{3}\right).
\end{equation}  
If we add the above two equations and use the relation
$$
\Pi_1(w)+\Pi_1\left(-\frac{1}{w}\right)=\frac{-\sqrt{3}\Gamma\left(\frac{1}{3}\right)^3}{\pi \sqrt[3]{2}},
$$
we get easily
\begin{equation}
F_c\left(Y(z)\right)+F_c\left(Y\left(-\frac{1}{z}\right)\right)=\frac{-\sqrt{3}\Gamma\left(\frac{1}{3}\right)^3}{\pi \sqrt[3]{2}}=:C_0.
\end{equation}
Solving this last equation with respect to $Y(-1/z)$ we get the result.\\
\\
\textbf{Examples.}\\
If $q=e^{2i\pi z}$, $Im(z),Im(w)>0$ and $Re(z),Re(w)\in\left(-\frac{1}{2},\frac{1}{2}\right)$, then with 
\begin{equation}
P(x)=2^{-1}x^{1/6}\sqrt{125+22x+x^2},
\end{equation}
we get (see [21])
\begin{equation}
Y(z)=\left(\frac{\eta(z)}{\eta(5z)}\right)^6=R\left(q\right)^{-5}-11-R\left(q\right)^5.
\end{equation}
Hence
\begin{equation}
2\int^{Y(w)}_{Y(z)}\frac{dt}{t^{1/6}\sqrt{125+22 t+t^2}}=
$$
$$
=-\left(\frac{2}{\sqrt[3]{4}}B\left(m^{*}(2z)^2;\frac{1}{6},\frac{2}{3}\right)-\frac{2}{\sqrt[3]{4}}B\left(m^{*}(2w)^2;\frac{1}{6},\frac{2}{3}\right)\right)
\end{equation}
2) For 
\begin{equation}
P(x)=\frac{\left(27+x^{12}\right)^{2/3}}{24x^5},
\end{equation}
we get after solving (116)
\begin{equation}
Y(z)=X\left(m^{*}(2z)\right)=\frac{\eta(z)}{\eta(3z)}.
\end{equation}
Hence
$$
12\int^{Y(w)}_{Y(z)}\frac{t^5}{(27+t^{12})^{2/3}}dt=
$$
\begin{equation}
=-\frac{2}{\sqrt[3]{4}}B\left(m^{*}(2w)^2;\frac{1}{6},\frac{2}{3}\right)+\frac{2}{\sqrt[3]{4}}B\left(m^{*}(2z)^2;\frac{1}{6},\frac{2}{3}\right)
\end{equation}
Hence it holds
\begin{equation}
\frac{2}{9}Y(z)^6\cdot {}_2F_1\left(\frac{1}{2},\frac{2}{3};\frac{3}{2};-\frac{Y(z)^{12}}{27}\right)+\frac{1}{\sqrt[3]{4}}B\left(m^{*}(2z)^2,\frac{1}{6},\frac{2}{3}\right)=-\frac{1}{2}C_0,
\end{equation}
for a constant $C$. Also
$$
\frac{2}{9}Y\left(-\frac{1}{z}\right)^6\cdot {}_2F_1\left(\frac{1}{2},\frac{2}{3};\frac{3}{2};-\frac{Y\left(-\frac{1}{z}\right)^{12}}{27}\right)+
$$
\begin{equation}
+\frac{2}{9}Y(z)^6\cdot {}_2F_1\left(\frac{1}{2},\frac{2}{3};\frac{3}{2};-\frac{Y(z)^{12}}{27}\right)=\frac{\sqrt{3}\Gamma\left(\frac{1}{3}\right)^3}{2\pi \sqrt[3]{2}}.
\end{equation}

\section{Construction of Hauptmodul functions of index   $N$}

In the same way as in above section, we define the function $X(x)$ as the solution of differential equation
\begin{equation}
X'(x)+2^{4/3}x^{-2/3}\left(1-x^2\right)^{-1/3}P\left(X(x)\right)=0.
\end{equation}  
Then the two solutions of the DE, ($N>0$):
\begin{equation}
Y'_{\mp}(z)=\pm 4\pi i \sqrt{N}\eta_D\left(\sqrt{N}z\right)^4P\left(Y_{\mp}(z)\right),
\end{equation}
are
\begin{equation}
Y_{-}(z)=X\left(m^{*}\left(2\sqrt{N}z\right)\right)
\end{equation}
and
\begin{equation} Y_{+}(z)=X\left(m^{*}\left(\frac{-2}{\sqrt{N}z}\right)\right)=Y_{-}\left(\frac{-1}{N z}\right),
\end{equation}
where $-\frac{1}{2\sqrt{N}}< Re(z)< \frac{1}{2\sqrt{N}}$ and $Im(z)>0$. Hence using Theorem 17, we get
\begin{equation}
\int^{Y_{\pm}(i\infty)}_{Y_{\pm}(z)}\frac{dt}{P(t)}=\mp\frac{2}{\sqrt[3]{4}}B\left(m^{*}\left(2z\sqrt{N}\right)^2;\frac{1}{6},\frac{2}{3}\right).
\end{equation}
Hence we have
\begin{equation}
\int^{Y_{+}(i\infty)}_{Y_{+}(z)}\frac{dt}{P(t)}+\int^{Y_{-}(i\infty)}_{Y_{-}(z)}\frac{dt}{P(t)}=0.
\end{equation}
Consequently we get
\begin{equation}
\int^{Y_{+}(i\infty)}_{Y_{-}\left(-1/(Nz)\right)}\frac{dt}{P(t)}+\int^{Y_{-}(i\infty)}_{Y_{-}(z)}\frac{dt}{P(t)}=0.
\end{equation}
If $F_c(x)=F(c,x)$ as in (133) and $c=Y_{+}(i\infty)$, $c'=Y_{-}(i\infty)$ we get
$$
F\left(c,Y_{-}\left(\frac{-1}{Nz}\right)\right)+F\left(c',Y_{-}\left(z\right)\right)=0\Leftrightarrow 
$$
$$
F\left(c,Y_{-}\left(\frac{-1}{Nz}\right)\right)+F\left(c,Y_{-}\left(z\right)\right)=F(c,c')\Leftrightarrow
$$
\begin{equation}
F_c\left(Y_{-}\left(\frac{-1}{Nz}\right)\right)+F_c\left(Y_{-}\left(z\right)\right)=C,
\end{equation}
where $C=F(c,c')$. But if we consider the ''regularized'' incomplete integral of $1/P(x)$ as 
\begin{equation}
F_{reg}(x)=\frac{1}{\int^{c'}_{c}\frac{dt}{P(t)}}\int^{x}_{c}\frac{dt}{P(t)},
\end{equation}
then
\begin{equation}
F_{reg}\left(Y_{-}\left(\frac{-1}{Nz}\right)\right)+F_{reg}\left(Y_{-}\left(z\right)\right)=1
\end{equation}
and
\begin{equation}
Y_{-}\left(\frac{-1}{Nz}\right)=F^{(-1)}_{reg}\left(1-F_{reg}\left(Y_{-}(z)\right)\right).
\end{equation}
\\

Assume $X(x)$ is solution of (149) and $Y_{\pm}(z)$ as in (151),(152), satisfying the DE (150). Then for every $z_1,z_2\in H$ such that $Y(z_1),Y(z_2)\in H$, we have from Theorem 21 below that
$$
\prod^{\infty}_{n=1}\left(\frac{1-e\left(nY(z_2)\right)}{1-e\left(nY(z_1)\right)}\right)^{-1/n\sum_{d|n}a_d\mu\left(n/d\right)}=\exp\left(2\pi i\int^{Y(z_2)}_{Y(z_1)}\left(\frac{1}{P(z)}-a_0\right)dz\right)\Leftrightarrow
$$

$$
\left(\frac{e\left(Y(z_2)\right)}{e\left(Y(z_1)\right)}\right)^{a_0}\prod^{\infty}_{n=1}\left(\frac{1-e\left(nY(z_2)\right)}{1-e\left(nY(z_1)\right)}\right)^{-1/n\sum_{d|n}a_d\mu\left(n/d\right)}=
$$
$$
=\exp\left(8\pi^2\int^{z_2\sqrt{N}}_{z_1\sqrt{N}}\eta\left(z\right)^4dz\right).
$$
Hence\\
\\
\textbf{Theorem 19.1}\\
If $X(x)$ is defined from
$$
X'(x)+2^{4/3}x^{-2/3}\left(1-x^2\right)^{-1/3}P\left(X(x)\right)=0,
$$
then we set $Y(z)=X\left(m^{*}\left(2z\sqrt{N}\right)\right)$. For such $Y$ holds 
$$
Y'(z)=\pm 4\pi i\sqrt{N}\eta\left(\sqrt{N}z\right)^4P\left(Y(z)\right).
$$
Also if $a_n$ are defined from
$$
\frac{1}{P(z)}=\sum^{\infty}_{n=0}a_nq^n\textrm{, }q=e(z)\textrm{, }Im(z)>0,\eqno{(159.1)}
$$
then one can see that exists constant $C$ such that
$$
e\left(a_0Y(z)\right)\prod^{\infty}_{n=1}\left(1-e\left(nY(z)\right)\right)^{-1/n\sum_{d|n}a_d\mu\left(n/d\right)}=
$$
$$=C\exp\left(8\pi^2\int^{z\sqrt{N}}_{i\infty}\eta\left(w\right)^4dw\right).\eqno{(159.2)}
$$
Also\small 
$$
e\left(a_0z\right)\prod^{\infty}_{n=1}\left(1-e\left(nz\right)\right)^{-1/n\sum_{d|n}a_d\mu\left(n/d\right)}=C\exp\left(8\pi^2\int^{Y^{(-1)}(z)\sqrt{N}}_{i\infty}\eta\left(w\right)^4dw\right).\eqno{(159.3)}
$$
\normalsize
\\
\textbf{Notes.}\\
We re-define $m(z)$ as 
$$2\pi i\int^{m(z)}_{i\infty}\eta(w)^4dw=z,\eqno{(159.4)}$$
then for every $z:Im(z)>0$, there exists integer $k$ such that\small
$$
m\left(\frac{i}{4\pi}\log\left(e\left(a_0Y(z)\right)\prod^{\infty}_{n=1}\left(1-e\left(nY(z)\right)\right)^{-1/n\sum_{d|n}a_d\mu\left(n/d\right)}\right)+\frac{k}{2}\right)=z\sqrt{N}
$$\normalsize
and if $Y(z)$ have inverse in the sense $Y(Y^{(-1)}(z))=z$, then\small
$$
Y^{(-1)}(z)=\frac{1}{\sqrt{N}}m\left(\frac{i}{4\pi}\log\left(e\left(a_0z\right)\prod^{\infty}_{n=1}\left(1-e\left(nz\right)\right)^{-1/n\sum_{d|n}a_d\mu\left(n/d\right)}\right)+\frac{k}{2}\right)\eqno{(159.5)}
$$\normalsize
\\
\textbf{Example.}\\
If $g(t)=\frac{1}{1+t^{\nu}}$ and $P(x)=g(e^{2\pi ix})$, then 
$$
X(x)=-6\sqrt[3]{2x}\cdot{}_2F_1\left(\frac{1}{6},\frac{1}{3};\frac{7}{6};x^2\right)
+\frac{i}{2\pi\nu}P_{L}\left(\frac{1}{2}e^{-12i\sqrt[3]{2x}\pi\nu \cdot{}_2F_1\left(\frac{1}{6},\frac{1}{3};\frac{7}{6};x^2\right)}\right).
$$
Hence $Y(z)=X\left(m^{*}\left(2z\right)\right)$. In this case the constant of integration have chosen such $Y'(z)=4\pi i\sqrt{N}\eta(\sqrt{N}z)^4P(Y(z))$.  The function $w=P_L(x)=W(x)$ is the Product logarithm and is defined as the solution of $we^w=x$.\\
Also 
$$
c_{\infty}=\lim_{h\rightarrow 0}X(h)=\frac{iW(1)}{2\pi\nu}
$$
and
$$
F(x)=\int^{x}_{c_{\infty}}\frac{dt}{P(t)}=-\frac{ie^{2\pi i \nu z}}{2\pi \nu}+z
$$
and
$$
F\left(Y\left(z\right)\right)+F\left(Y\left(-\frac{1}{z}\right)\right)=\frac{\sqrt{3}\Gamma\left(\frac{1}{3}\right)^3}{\sqrt[3]{2}\pi}.
$$
Also we have 
$$
e\left(Y(z)\right)\prod_{\scriptsize\begin{array}{cc}
n\geq 1\\
n\equiv0(\nu)
\end{array}\normalsize}\left(1-e(nY(z))\right)^{-\mu(n/\nu)/n}=\exp\left(8\pi^2\int^{z}_{i\infty}\eta\left(w\right)^4dw\right).
$$
\\
\textbf{Example.}\\
Assume $\frac{1}{P(z)}=\sum^{\infty}_{n=1}\chi(\sqrt[\nu]{n})X_{\nu}(n)nq^{n}=\sum^{\infty}_{n=1}\chi(n)n^{\nu}q^{n^{\nu}}$. Then
$$
\prod^{\infty}_{n=1}\left(1-e\left(nY_{N}(z)\right)\right)^{-1/n\sum_{d^{\nu}|n}\chi(d)d^{\nu}\mu\left(n/d^{\nu}\right)}=
$$
$$
=e^C\cdot\exp\left(8\pi^2\int^{z\sqrt{N}}_{i\infty}\eta\left(w\right)^4dw\right).
$$
Hence there is a function $Y_{N}(z)$ such that
$$
\exp\left(\sum^{\infty}_{n=1}\chi(n)e\left(n^{\nu}Y_{N}(z)\right)\right)=e^C\cdot\exp\left(8\pi^2\int^{z\sqrt{N}}_{i\infty}\eta\left(w\right)^4dw\right)\Leftrightarrow
$$
$$
\sum^{\infty}_{n=1}\chi(n)e\left(n^{\nu}Y_{N}(z)\right)=8\pi^2\int^{z\sqrt{N}}_{i\infty}\eta(w)^4dw+C+2\pi i k\textrm{, }k\in\textbf{Z}.
$$
Or ''equivalently''
$$
Y_{N}(z)=\theta^{{\{\nu\}}{(-1)}}_{\chi}\left(8\pi^2\int^{z\sqrt{N}}_{i\infty}\eta(w)^4dw+C'\right).
$$
Also $Y_N(z)$ satisfies the modular relation
$$
\theta^{\{\nu\}}_{\chi}\left(Y_{N}\left(\frac{-1}{N z}\right)\right)+\theta^{\{\nu\}}_{\chi}\left(Y_{N}\left(z\right)\right)=C_1,
$$
where $C_1$ is constant.\\
\\
\textbf{Example.}\\
If $q=e(z)$, $-\frac{1}{2\sqrt{N}}<Re(z)<\frac{1}{2\sqrt{N}}$ and $Im(z)>0$, then for 
$$
P(x)=\frac{\left(27+x^{12}\right)^{2/3}}{12x^5},
$$
we get after solving (147)
$$
Y(z)=X\left(m^{*}\left(2z\right)\right)=\frac{\eta_D(z)}{\eta_D(3z)}.
$$
Hence if
\begin{equation}
Y_1(z)=X\left(m^{*}\left(2\sqrt{N}z\right)\right)=\frac{\eta_D\left(\sqrt{N} z\right)}{\eta_D\left(3\sqrt{N}z\right)},
\end{equation}
then it holds
$$
Y_1\left(-\frac{1}{Nz}\right)^6\cdot {}_2F_1\left(\frac{1}{2},\frac{2}{3};\frac{3}{2};-\frac{Y_1\left(-\frac{1}{Nz}\right)^{12}}{27}\right)+
$$
\begin{equation}
+Y_1(z)^6\cdot {}_2F_1\left(\frac{1}{2},\frac{2}{3};\frac{3}{2};-\frac{Y_1(z)^{12}}{27}\right)=\frac{9\sqrt{3}\Gamma\left(\frac{1}{3}\right)^3}{4\pi \sqrt[3]{2}}.
\end{equation}
\\
\textbf{Example.}\\
Assume $a_n=n$, then $1/P(q)=\sum^{\infty}_{n=0}a_nq^n=\frac{(q-1)^2}{q}$. Hence $P_1(t)=P\left(e^{2\pi i t}\right)=-4\sin(\pi t)^2$ and equation
$$
X'(x)+2^{4/3}x^{-2/3}(1-x^2)^{-1/3}P_1(X(x))=0,
$$
have solution
$$
X(x)=-\frac{1}{\pi}\textrm{arccot}\left(24\pi \sqrt[3]{2x}\cdot {}_2F_1\left(\frac{1}{6},\frac{1}{3};\frac{7}{6};x^2\right)\right).
$$
Also $\lim_{h\rightarrow +\infty}Y(ih)=\lim_{h\rightarrow 0}X(h)=-\frac{1}{2}$ and
$$
\int^{x}_{-1/2}\frac{dt}{P_1(t)}=\frac{\cot(\pi x)}{4\pi}\textrm{, }Re(x)\leq 0\textrm{, }Im(x)\neq 0.
$$
Hence if $Y(z)=X\left(m^{*}(2\sqrt{N}z)\right)$, then
$$
Y\left(-\frac{1}{Nz}\right)=-\pi^{-1}\textrm{arccot}\left(\cot\left(\pi Y(z)\right)+4\pi C_0\right)\textrm{, }C_0=\frac{\sqrt{3}\Gamma\left(\frac{1}{3}\right)^3}{\sqrt[3]{2}\pi}.
$$
Also $\sum_{d|n}a_d\mu(n/d)=\phi(n)$ and exists constant $C$ such that
$$
\prod^{\infty}_{n=1}\left(1-e\left(nY(z)\right)\right)^{-\phi(n)/n}=C\cdot \exp\left(8\pi^2\int^{z\sqrt{N}}_{i\infty}\eta(w)^4dw\right).
$$

\section{A more traditional way to study $Q_{\{a,p\}}(x)$}

\textbf{Proposition 1.}\\
If $x$ is positive real number and $f$ is analytic in $(-1,1)$ with $f(0)=0$, then
\begin{equation}
\exp\left(\int^{x}_{+\infty}f(e^{-t})dt\right)=\prod^{\infty}_{n=1}(1-e^{-nx})^{\frac{1}{n}\sum_{d|n}\frac{f^{(d)}(0)}{d!}\mu\left(\frac{n}{d}\right)},
\end{equation}
where $\mu$ is the Moebius-$\mu$ arithmetic function (see [17]) and take the values $(-1)^r$ when $n$ square free and product of $r$ primes, else is $0$. Also $\mu(1)=1$.\\
\\
\textbf{Proof.}\\
Because $f(0)=0$ and $f$ analytic in $(-1,1)$, the integral $\int^{x}_{+\infty}f(e^{-t})dt$ exists for every $x>0$. We assume that exists arithmetic function $X(n)$ such that: 
\begin{equation}
\exp\left(\int^{x}_{+\infty}f(e^{-t})dt\right)=\prod^{\infty}_{n=1}(1-e^{-nx})^{X(n)}.
\end{equation}
We will determinate this function $X$.\\
Taking logarithms in both sides of (163) we have
$$
\int^{x}_{+\infty}f(e^{-t})dt=\sum^{\infty}_{n=1}X(n)\log(1-e^{-nx})=-\sum^{\infty}_{n=1}X(n)\sum^{\infty}_{m=1}\frac{e^{-mnx}}{m}=
$$
$$
=-\sum^{\infty}_{n,m=1}X(n)n\frac{e^{-mnx}}{mn}=-\sum^{\infty}_{n=1}\frac{e^{-nx}}{n}\sum_{d|n}X(d)d.\eqno{:(A)}$$
Derivating (A) we get
$$
f(x)=\sum^{\infty}_{n=1}e^{-nx}\sum_{d|n}X(d)d.\eqno{:(B)}
$$
But from analytic property of $f$ in $(-1,1)$ we have 
$$
f(x)=\sum^{\infty}_{n=1}\frac{f^{(n)}(0)}{n!}x^n
$$ 
and consequently
$$
f(e^{-x})=\sum^{\infty}_{n=1}\frac{f^{(n)}(0)}{n!}e^{-nx}.
$$
Therefore from (B) and the above relation it must be 
$$
\frac{f^{(n)}(0)}{n!}=\sum_{d|n}X(d)d.
$$
By applying the Moebius inversion theorem (see [17]) we get
$$
X(n)=\frac{1}{n}\sum_{d|n}\frac{f^{(d)}(0)}{d!}\mu\left(\frac{n}{d}\right).
$$
This completes the proof.\\ 
\\
\textbf{Theorem 20.}\\
Let $|q|<1$, then 
\begin{equation}
e^{-f(q)}=\prod^{\infty}_{n=1}\left(1-q^n\right)^{\frac{1}{n}\sum_{d|n}\frac{f^{(d)}(0)}{\Gamma(d)}\mu\left(\frac{n}{d}\right)}.
\end{equation}
\\
\textbf{Proof.}\\
Setting where $\frac{f^{(n)}(0)}{n!}=\frac{f_1^{(n)}(0)}{n!}n$ and using Proposition 1, we get imediately the result.\\
\\
\textbf{Theorem 21.}\\
If $z_1,z_2$ are complex numbers in $\textbf{H}$ (the upper half plane) and 
\begin{equation}
f(z)=\sum^{\infty}_{n=1}a_nq^n\textrm{, }q=e(z)\textrm{, }Im(z)>0,
\end{equation}
holomorphic also in $\textbf{H}$, then
\begin{equation}
\exp\left(2\pi i \int^{z_2}_{z_1}f(z)dz\right)=\prod^{\infty}_{n=1}\left(\frac{1-q_2^n}{1-q_1^n}\right)^{-X(n)}, 
\end{equation} 
where $q_j=e(z_j)$, $j=1,2$ and 
\begin{equation}
X(n)=\frac{1}{n}\sum_{d|n}a_d\mu\left(\frac{n}{d}\right).
\end{equation}
\\
\textbf{Corollary.}\\
Assume that exists a function $f(z)$ and constants $k,N,\epsilon$ such that for all $z\in\textbf{H}$ hold
\begin{equation}
\exp\left(2\pi i \int^{-1/(Nz)}_{z}f\left(w\right)dw\right)=\epsilon z^k.
\end{equation}
Further if
\begin{equation}
\int^{z+1}_{z}f\left(w\right)dw=0,
\end{equation}
then 
\begin{equation}
\exp\left(2\pi i\int^{z}_{c_0}f(w)dw\right)
\end{equation}
is a modular form of weight $k$ in $\Gamma(N)$. Also (from Theorem 21) the function
\begin{equation}
\phi(z)=\prod^{\infty}_{n=1}(1-q^n)^{-X(n)},
\end{equation}  
is a modular form of weight $k$ in $\Gamma(N)$.\\
\\

Taking the logarithms and derivating both sides of (168), we can write 
$$
f\left(\frac{-1}{Nz}\right)\frac{1}{Nz^2}-f\left(z\right)=\frac{k}{2\pi i z}.
$$
Hence
$$
f\left(\frac{-1}{Nz}\right)\frac{-1}{N z}+zf\left(z\right)=-\frac{k}{2\pi i}.
$$
If we set 
\begin{equation}
g(z)=-\frac{2\pi i}{k}zf\left(z\right),
\end{equation}
then
\begin{equation}
g\left(\frac{-1}{Nz}\right)+g(z)=1.
\end{equation}
We can write
\begin{equation}
g\left(-\frac{1}{\sqrt{N}z}\right)+g\left(\frac{z}{\sqrt{N}}\right)=1
\end{equation}
and if
\begin{equation}
h(z)=g\left(\frac{z}{\sqrt{N}}\right),
\end{equation}
then
\begin{equation}
h\left(-\frac{1}{z}\right)+h(z)=1.
\end{equation}
Hence we get the next corollary.\\
\\
\textbf{Corollary.}\\
Let $h(z)$ be a function such that for all $z\in\textbf{H}$ we have
\begin{equation}
h\left(\frac{-1}{z}\right)+h(z)=1\textrm{,  }\frac{h\left(z+\sqrt{N}\right)}{z+\sqrt{N}}=\frac{h\left(z\right)}{z}.
\end{equation}
Then the function
\begin{equation}
f(z)=-\frac{k}{2\pi i z}h\left(z\sqrt{N}\right)
\end{equation}
have Fourier expansion
\begin{equation}
f(z)=\sum^{\infty}_{n=0}a_nq^n\textrm{, }q=e(z)\textrm{, }Im(z)>0,
\end{equation}
and if
\begin{equation}
X(n)=\frac{1}{n}\sum_{d|n}a_d\mu\left(\frac{n}{d}\right),
\end{equation}
the function 
\begin{equation}
\phi(z)=\prod^{\infty}_{n=1}\left(1-q^n\right)^{-X(n)},
\end{equation}
is a modular form of weight $k$ in a certain group $\Gamma(N)$. Also holds the next representation
\begin{equation}
\phi(z)=\exp\left(2\pi i\int^{z}_{i\infty}f(w)dw\right).
\end{equation}
\\
\textbf{Example 1.}\\
Assume that
\begin{equation}
E_4(z)=q^A\prod^{\infty}_{n=1}\left(1-q^n\right)^{-X(n)}.
\end{equation}
For to evaluate $A$ and $X(n)$, we write $M(q)=E_4(z)=\exp\left(f(q)\right)$, then from
$$
q\frac{dM}{dq}=\frac{LM-N}{3}
$$
and
$$
\frac{J'(q)}{J(q)}=-\frac{N}{qM},
$$
we get
$$
\frac{M'(q)}{M(q)}=\frac{L}{3q}-\frac{N}{3qM}=\frac{L}{3q}+\frac{J'(q)}{3J(q)},
$$
where $j(z)=J(q)$, $q=e(z)$, $Im(z)>0$. Hence writing 
$$
\frac{M'(q)}{M(q)}=\sum^{\infty}_{n=1}A_{n-1}q^n,
$$
we have
\begin{equation}
A_{n}=\frac{1}{3}c_n-8\sigma_1(n+1),
\end{equation}
where $c_n$ are the series coefficents of $J'(q)/J(q)$. Hence
$$
M(q)=\exp\left(\sum^{\infty}_{n=1}A_{n-1}\frac{q^n}{n}\right)=\prod^{\infty}_{n=1}\left(1-q^n\right)^{-X(n)},
$$
where
$$
-X(n)=\frac{1}{n}\sum_{d|n}A_{d-1}\mu\left(\frac{n}{d}\right)=8-\frac{1}{3n}\sum_{d|n}c_{d-1}\mu\left(\frac{n}{d}\right).
$$
Hence
\begin{equation}
E_4(z)=\prod^{\infty}_{n=1}\left(1-q^n\right)^{8-\frac{1}{3n}\sum_{d|n}c_{d-1}\mu(n/d)}.
\end{equation}
\\
\textbf{Example 2.}\\
Assume that $X(n)=1$ and $A=1/24$, then $\phi(z)=\eta(z)$, where $\eta(z)=q^{1/24}\prod^{\infty}_{n=1}\left(1-q^n\right)$ is the Dedekind eta function. This function have modular properties. i.e If $ad-bc=1$, then exist $\epsilon=\epsilon(a,b,c,d)$ and $\epsilon^{24}=1$ such that 
$$
\eta(\sigma(z))=\epsilon(a,b,c,d)(cz+d)^{1/2}\eta(z)\textrm{, }\forall z\in\textbf{H}.
$$
Hence if we assume the function
$$
f(z)=\sum^{\infty}_{n=1}\sigma_1(n)q^n\textrm{, }q=e(z)\textrm{, }Im(z)>0,
$$
then
$$
X(n)=1=\frac{1}{n}\sum_{d|n}\sigma_1(n)\mu(n/d)
$$
and $\psi(z)=\exp\left[2\pi i(z/24-F(z))\right]$, where $F'(z)=f(z)$, behaves exactly as $\eta(z)$ i.e. is a modular form of weight 1/2 and  
$$
\psi(\sigma(z))=\epsilon(a,b,c,d)(cz+d)^{1/2}\psi(z)\textrm{, }\forall z\in\textbf{H}.
$$
Actualy it is
$$
\eta(z)=\exp\left[2\pi i\left(z/24-F(z)\right)\right].
$$
Hence in better detail
\begin{equation}
\exp\left(2\pi i\int^{z_2}_{z_1}E_2(z)dz\right)=\frac{\Delta(z_2)}{\Delta(z_1)},
\end{equation}
where $\Delta(z)=\eta_{D}(z)^{24}$.\\
More generaly if $\nu$ is even positive integer and   
$$
F_{\nu}(z)=-\frac{1}{2\pi i}\frac{2\nu}{B_{\nu}}\sum^{\infty}_{n=1}\sigma_{\nu-1}(n)\frac{q^n}{n},
$$
then
$$
z_2+F_{2\nu}(z_2)-\left(z_1+F_{2\nu}(z_1)\right)=\int^{z_2}_{z_1}E_{2\nu}(z)dz
$$
and
\begin{equation}
\frac{q_2}{q_1}\prod^{\infty}_{n=1}\left(\frac{1-q_2^n}{1-q_1^n}\right)^{4n^{2\nu-2}\nu/B_{2\nu}}=\exp\left(2\pi i\int^{z_2}_{z_1}E_{2\nu}(z)dz\right),
\end{equation}
since 
$$
\frac{1}{n}\sum_{d|n}\sigma_{\nu-1}(d)\mu\left(\frac{n}{d}\right)=n^{\nu-2}.
$$
\\
\textbf{Example 3.}\\
If $\lambda(n)$ is Liouville's lambda arithmetical function and if $q_j=e(z_j)$, $Im(z_j)>0$, $j=1,2$, then
\begin{equation}
\exp\left(2\pi i\int^{z_2}_{z_1}\theta_3(z)dz\right)=\frac{q_2}{q_1}\prod^{\infty}_{n=1}\left(\frac{1-q_2^n}{1-q_1^n}\right)^{-2\lambda(n)/n},
\end{equation}
where
\begin{equation}
\theta_3(z)=\sum^{\infty}_{n=-\infty}q^{n^2}\textrm{, }q=e(z)\textrm{, }Im(z)>0.
\end{equation}
\\
\textbf{Proof.}\\
Use Theorem 21 and the identities
$$
\frac{\theta_{3}(q)-1}{2}=\sum^{\infty}_{n=1}X_2(n)q^n=\sum^{\infty}_{n=1}q^{n^2}\textrm{, }|q|<1,
$$
$$
\lambda(n)=\sum_{d^2|n}\mu\left(\frac{n}{d^2}\right),
$$
where $\mu(n)$ is the Moebious $\mu$ arithmetical function.\\
\\

More generaly one can see that\\
\\
\textbf{Theorem 22.}\\
If $q=e(z)$, $Im(z)>0$ and define the next generalization of theta function as
\begin{equation}
\psi_{\nu}(z):=\sum^{\infty}_{n=1}q^{n^{\nu}}\textrm{, }\nu=2,3,4,\ldots,
\end{equation}
then holds
\begin{equation}
\exp\left(2\pi i\int^{z}_{i\infty}\psi_{\nu}(w)dw\right)=\prod^{\infty}_{n=1}\left(1-q^n\right)^{-\lambda_{\nu}(n)/n},
\end{equation}
where 
\begin{equation}
\lambda_{\nu}(n):=\sum_{d^{\nu}|n}\mu\left(\frac{n}{d^{\nu}}\right),
\end{equation}
is the generalized Liouville function.\\
For this function also holds
\begin{equation}
\sum_{d|n}\lambda_{\nu}(d)=X_{\nu}(n):=\left\{
\begin{array}{cc}
	1\textrm{, if }\exists m\in\textbf{N}:n=m^{\nu}\\
	0\textrm{, else }
\end{array}
\right\}.
\end{equation}
Also if $(n,m)=1$, then
\begin{equation}
\lambda_{\nu}(nm)=\lambda_{\nu}(n)\lambda_{\nu}(m).
\end{equation}
\begin{equation}
\lambda_{\nu}\left(n^{\nu}\right)=1.
\end{equation}
\\
\textbf{Thorem 23.}\\
If $|q|<1$, then
\begin{equation}
\psi_{\nu}(z)=\sum^{\infty}_{n=1}q^{n^{\nu}}=\sum^{\infty}_{n=1}X_{\nu}(n)q^n\textrm{, }|q|<1.
\end{equation}
Then
\begin{equation}
X_{\nu}(nm)X_{\nu}\left(\textrm{gcd}(n,m)\right)=X_{\nu}(n)X_{\nu}(m)\textrm{, }\forall n,m\in\{1,2,\ldots\}
\end{equation}
and $X_{\nu}(n)$ have Dirichlet series
\begin{equation}
L\left(X_{\nu},s\right)=\sum^{\infty}_{n=1}\frac{X_{\nu}(n)}{n^s}=\zeta(\nu s),
\end{equation}
where $\zeta(s)$ is the Riemann's zeta function.  Also
\begin{equation}
\exp\left(\psi_{\nu}(z)\right)=\prod^{\infty}_{n=1}\left(1-q^n\right)^{-1/n\sum_{d^{\nu}|n}d^{\nu}\mu\left(n/d^{\nu}\right)}
\end{equation}
\\
\textbf{Proof.}\\
The proofs of (197) and (198) are easy. For to prove (199), we have
$$
\psi_{\nu}'(z)=2\pi i\sum^{\infty}_{n=1}X_{\nu}(n)nq^{n}.
$$
Hence from Theorem 21 we have
$$
\exp\left(2\pi i \int^{z}_{i\infty}\psi_{\nu}'(w)dw\right)=\prod^{\infty}_{n=1}\left(1-q^n\right)^{-2\pi i/n \sum_{d|n}X_{\nu}(d)d\mu(n/d)}.
$$
From this last relation we get (199).\\
\\
\textbf{Theorem 24.}\\
Suppose that $\nu=2,3,4,\ldots$ and $q=e(z)$, $Im(z)>0$. Then for any analytic function $g(z)$, $z\in\textbf{C}$, such that 
\begin{equation}
\left|g^{(k)}(0)\right|\leq C A^k\textrm{, }\forall k=0,1,2,\ldots
\end{equation}
and $A>0$, we have
\begin{equation}
\exp\left(\sum^{\infty}_{n=1}g\left(2\pi i n^{\nu}\right)q^{n^{\nu}}\right)=\prod^{\infty}_{n=1}\left(1-q^n\right)^{-X_g(n)}\textrm{, }\forall q:|q|<e^{-2\pi A},
\end{equation}
where 
\begin{equation}
X_g(n)=\frac{1}{n}\sum_{d^{\nu}|n}g\left(2\pi i d^{\nu}\right)d^{\nu}\mu\left(\frac{n}{d^{\nu}}\right).
\end{equation}
\\
\textbf{Proof.}\\
Since (200) holds we can write
$$
|g(z)|=\left|\sum^{\infty}_{k=0}\frac{g^{(k)}(0)}{k!}z^k\right|\leq \sum^{\infty}_{k=0}\frac{\left|g^{(k)}(0)\right|}{k!}|z|^k\leq C\sum^{\infty}_{k=0}\frac{A^k|z|^k}{k!}= Ce^{A|z|}.
$$
Hence
$$
\left|g(z)\right|\leq C e^{A|z|}\textrm{, }\forall z\in\textbf{C}.\eqno({a})
$$
As in Theorem 23, we take the $(k+1)-$th derivative of $\psi_{\nu}(z)$ with respect to $z$. We have
$$
\psi_{\nu}^{(k+1)}(z)=\sum^{\infty}_{n=1}(2\pi i n^{\nu})^{k+1}q^{n^{\nu}}=\sum^{\infty}_{n=1}(2\pi i n)^{k+1}X_{\nu}(n)q^{n}.
$$
Hence from Theorem 21, we have
$$
2 \pi i\int^{z}_{i\infty}\psi^{(k+1)}_{\nu}(w)dw=-\sum^{\infty}_{n=1}\frac{1}{n}\sum_{d^{\nu}|n}(2\pi i d^{\nu})^{k+1}\mu\left(\frac{n}{d^{\nu}}\right)\log(1-q^n).
$$
Hence
$$
2 \pi i\psi^{(k)}_{\nu}(z)=-2\pi i\sum^{\infty}_{n=1}\frac{1}{n}\sum_{d^{\nu}|n}(2\pi i d^{\nu})^{k}d^{\nu}\mu\left(\frac{n}{d^{\nu}}\right)\log(1-q^n)=
$$
$$
=-2\pi i\sum^{\infty}_{n=1}\frac{1}{n}\sum_{d|n}(2\pi i d)^{k}X_{\nu}(d) d\mu\left(\frac{n}{d}\right)\log(1-q^n),
$$
for all $k=0,1,2,\ldots$. Note that
$$
\left|\frac{1}{n}\sum_{d^{\nu}|n}(2\pi i d^{\nu})^{k}d^{\nu}\mu\left(\frac{n}{d^{\nu}}\right)\right|= \left|\frac{(2\pi)^k}{n}\sum_{d|n}X_{\nu}(d)d^{k+1}\mu\left(\frac{n}{d}\right)\right|\leq
$$
$$
\leq\frac{(2\pi)^k}{n}\sum_{d|n}d^{k+1}=\frac{(2\pi)^k\sigma_{k+1}(n)}{n}\leq C (2\pi n)^k n^{\epsilon}\textrm{, }\forall\epsilon>0\eqno{(b)}
$$ 
and 
$$
\left|\log(1-q^n)\right|\leq \frac{|q|^n}{1-|q|^n}\textrm{, }|q|<1.\eqno{(c)}
$$
Also  
$$
M_{n,k}=\left|\frac{g^{(k)}(0)}{k!}\left(2\pi i n\right)^{k}X_{\nu}(n)q^{n}\right|\leq C\frac{(2\pi A n)^k}{k!}|q|^{n}
$$
and
$$
\sum^{\infty}_{n=1}\sum^{\infty}_{k=0}M_{n,k}\leq C\sum^{\infty}_{n=1}\sum^{\infty}_{k=0}\frac{(2\pi A n)^k}{k!}|q|^{n}=C\sum^{\infty}_{n=1}e^{2\pi A n}|q|^{n}=
$$
$$
=C\sum^{\infty}_{n=1}\left(e^{2\pi A }|q|\right)^{n}<\infty\textrm{, when }|q|<e^{-2\pi A}.\eqno{(d)}
$$
Hence
$$
\sum^{\infty}_{n=1}g\left(2\pi i n^{\nu}\right)q^{n^{\nu}}=\sum^{\infty}_{n=1}g\left(2\pi i n\right)X_{\nu}(n)q^{n}=
$$
$$
=\sum^{\infty}_{n=1}\sum^{\infty}_{k=0}\frac{g^{(k)}(0)}{k!}\left(2\pi i n\right)^{k}X_{\nu}(n)q^{n}=\sum^{\infty}_{k=0}\frac{g^{(k)}(0)}{k!}\psi^{(k)}_{\nu}(z)
$$
and
$$
\psi^{(k)}_{\nu}(z)=-\sum^{\infty}_{n=1}\frac{1}{n}\sum_{d^{\nu}|n}(2\pi i d^{\nu})^{k}d^{\nu}\mu\left(\frac{n}{d^{\nu}}\right)\log(1-q^n).
$$
Hence easily we conclude that if $|q|<e^{-2\pi A}$, then
$$
\sum^{\infty}_{n=1}g\left(2\pi i n^{\nu}\right)q^{n^{\nu}}=-\sum^{\infty}_{n=1}\frac{1}{n}\sum_{d^{\nu}|n}g\left(2\pi i d^{\nu}\right)d^{\nu}\mu\left(\frac{n}{d^{\nu}}\right)\log(1-q^n).
$$
\\
\textbf{Notes.}\\
If 
$$
X_g(n)=\frac{1}{n}\sum_{d^{\nu}|n}g\left(2\pi i d^{\nu}\right)d^{\nu}\mu\left(\frac{n}{d^{\nu}}\right),
$$
then
$$
X_g(n)=\frac{1}{n}\sum_{d|n}g(2\pi i d)X_{\nu}(d)d\mu\left(\frac{n}{d}\right)\Leftrightarrow\sum_{d|n}X_g(d)d=g(2\pi i n)nX_{\nu}(n).
$$
\\
\textbf{Theorem 25.}\\
If $q=e(z)$, $Im(z)>0$ and $\chi(n)$ is any arithmetical function such that 
$$
\sum^{\infty}_{n=1}\left|\chi(n)\right|\cdot|q|^{n^{\nu}}<\infty,
$$ 
then
\begin{equation}
\exp\left(\sum^{\infty}_{n=1}\chi(n)q^{n^{\nu}}\right)=\prod^{\infty}_{n=1}\left(1-q^n\right)^{-X_{\chi}(n)},
\end{equation}
where
\begin{equation}
X_{\chi}(n)=\frac{1}{n}\sum_{d^{\nu}|n}\chi(d)d^{\nu}\mu\left(\frac{n}{d^{\nu}}\right).
\end{equation}
\\
\textbf{Proof.}\\
If we set 
$$
f(q)=\sum^{\infty}_{n=1}a_nq^n,
$$
where $a_n=\frac{f^{(n)}(0)}{n!}=\chi(n)X_{\nu}(n)$ in (164) of Theorem 20 and then $\chi(n^{\nu})\rightarrow\chi(n)$, we get the result.\\ 
\\
\textbf{Theorem 26.}\\
If $q=e(z)$, $Im(z)>0$ and $\nu=2,3,\ldots$, then
\begin{equation}
\psi_{\nu}(z)=\sum^{\infty}_{n=1}q^{n^{\nu}}=\sum^{\infty}_{n=1}\frac{\lambda_{\nu}(n)q^n}{1-q^n}.
\end{equation}
\\
\textbf{Proof.}\\
$$
\psi_{\nu}(z)=\sum^{\infty}_{n=1}X_{\nu}(n)q^n=\sum^{\infty}_{n=1}q^n\sum_{d|n}\lambda_{\nu}(d)=
$$
$$
=\sum^{\infty}_{n,m=1}q^{nm}\lambda_{\nu}(n)=\sum^{\infty}_{n=1}\frac{\lambda_{\nu}(n)q^n}{1-q^n}.
$$
\\
\textbf{Proposition 2.}\\
If $q=e^{-2x}$, $x>0$ and $k,h$ are integers such $k>0$ and $k>h$ then
\begin{equation}
\log\left(\sum^{\infty}_{n=-\infty}(-1)^nq^{kn^2+hn}\right)=-\sum^{\infty}_{n=1}\left(\sum_{d|n}\chi_{k,h}(d)d\right)\frac{q^n}{n},
\end{equation}
where $\chi_{k,h}(n)$ is that of (208) below.\\
\\
\textbf{Proof.}\\
Assume the Jacobi's triple product identity (see [16] pg.169-172 and Exercise 3 pg.178)
\begin{equation}
\sum^{\infty}_{n=-\infty}(-1)^nq^{kn^2+hn}=\prod_{n=0}^{\infty}\left(1-q^{2kn+k-h}\right)\left(1-q^{2kn+k+h}\right)\left(1-q^{2kn+2k}\right),
\end{equation}
where $|q|<1$, $k>0$.\\Setting
\begin{equation} 
\chi_{k,h}(n):=\left\{
\begin{array}{cc}
  
	1 \textrm{ if } n\equiv 
	 
	0,k+h,k-h(mod2k)
\\	
	
	0 \textrm{ otherwise }	 

\end{array}
\right\}, 
\end{equation}
we can write
\begin{equation}
\sum^{\infty}_{n=-\infty}(-1)^n q^{kn^2+hn}=\prod^{\infty}_{n=1}\left(1-q^{n}\right)^{\chi_{k,h}(n)}
\end{equation}
Recall now Lemma 1 and Moebius inversion formula to write
$$
\exp\left(-\sum^{\infty}_{n=1}\frac{q^n}{n}\sum_{d|n}\chi_{k,h}(d)d\right)=\sum^{\infty}_{n=-\infty}(-1)^nq^{kn^2+hn}
$$
and hence (206).\\
\\
\textbf{Theorem 27.}\\
Let $q=e^{-\pi\sqrt{r}}$, $r>0$ and $a,p$ integers such that $p\geq a$ and $p>0$. If the theta function
\begin{equation}
\sum^{\infty}_{n=-\infty}(-1)^n q^{p n^2/2+(p-2a) n/2}=q^{A}\eta\left(q^{p}\right)Q_{\{a,p\}}\left(k_r\right)
\end{equation}
have algebraic part $Q_{\{a,p\}}(k_r)$, then
$$
\log\left(Q_{\{a,p\}}\left(k_{r}\right)\right)=
-A\log q-\sum^{\infty}_{n=1}\left(\sum_{\scriptsize 
\begin{array}{cc}
	d|n\\
	d\equiv\pm a(p)
\end{array}
\normalsize}d\right)\frac{q^{n}}{n}
$$
\begin{equation}
=-A\log q-\sum^{\infty}_{n=1}\left(\sum_{\scriptsize 
\begin{array}{cc}
	d|n\\
	d\equiv\pm a(p)
\end{array}
\normalsize}\frac{1}{d}\right)q^{n},
\end{equation}
where
\begin{equation}
A=-\frac{p}{12}+\frac{a}{2}-\frac{a^2}{2p}.
\end{equation} 
\\
\textbf{Proof.}\\
Let $X_{\{a,p\}}(n)$ be such that
\begin{equation} 
X_{\{a,p\}}(n):=\left\{
\begin{array}{cc}
  
	1 \textrm{ if } n\equiv 
	 
	0,a,p-a(\textrm{mod}p)
\\	
	
	0 \textrm{ otherwise }	 

\end{array}
\right\}. 
\end{equation}
From relation (30) we get 
\begin{equation}
\log\left(Q_{\{a,p\}}\left(k_{r}\right)\right)=-A\log q+\log\left(\sum^{\infty}_{n=-\infty}(-1)^nq^{pn^2/2+(p-2a)n/2}\right)-\log\left(\eta\left(q^{p}\right)\right).
\end{equation}
From Proposition 2 we have
\begin{equation}
\log\left(\sum^{\infty}_{n=-\infty}(-1)^nq^{pn^2/2+(p-2a)n/2}\right)=-\sum^{\infty}_{n=1}\left(\sum_{d|n}X_{\{a,p\}}(d)d\right)\frac{q^{n}}{n}.
\end{equation}
Also there holds
\begin{equation}
\log\left(\eta(q)\right)=\sum^{\infty}_{n=1}\log\left(1-q^n\right)=-\sum^{\infty}_{n,m=1}\frac{q^{nm}}{m}=-\sum^{\infty}_{n=1}\sigma_{-1}(n)q^n.
\end{equation} 
From (214),(215) and (216) we get the result. $qed$\\
\\

Using complex $q=e(z)$, $Im(z)>0$ and working as with $\theta_4(a,b;q)$, we get
\begin{equation}
\vartheta_3\left(\pi zt,e(az)\right)=q^{a/12-t^2/(4a)}\eta\left(q^{2a}\right)Q^{\{3\}}_{\{a,t\}}\left(m(q)\right),
\end{equation}
where
\begin{equation}
Q^{\{3\}}_{\{a,t\}}(m(q))=q^{-a/12+t^2/(4a)}\exp\left(-\sum^{\infty}_{n=1}q^n\sum_{\scriptsize
\begin{array}{cc}
AB=n\\
B\equiv \pm (a-t)(mod 2a)	
\end{array}\normalsize}\frac{(-1)^A}{A}\right).
\end{equation} 
But holds the following modular identity (see [2]):
\begin{equation}
\vartheta_3\left(\pi t' z',e(a'z')\right)=\sqrt{-2iaz}\exp\left(\frac{i\pi t^2 z}{2a}\right)\vartheta_3\left(\pi t z,e(az)\right),
\end{equation}
where
\begin{equation}
a'=1/a\textrm{, }z'=-1/(4z)\textrm{, }t'=2tz/a.
\end{equation}
Hence in general for the function $F_3(a,t;z):=Q^{\{3\}}_{\{a,t\}}(m(q))$, $q=e(z)$ holds
$$
\frac{F_3(a',t';z')}{F_3(a,t;z)}=\sqrt{-2iaz}\exp\left(\frac{i\pi t^2 z}{2a}\right)\frac{\eta_D(2az)}{\eta_D\left(\frac{-1}{2az}\right)},
$$
where $\eta_D(z)$, $Im(z)>0$ is the Dedekind's eta function. Using the next functional equation: 
\begin{equation}
\eta_D\left(-1/z\right)=\sqrt{-iz}\cdot\eta_D(z),
\end{equation}
we finally arrive to\\
\\
\textbf{Theorem 28.}\\
Let $a>0$, $q=e(z)$, $Im(z)>0$ and
\begin{equation}
\vartheta_3\left(\pi t z,e(az)\right)=q^{a/12-t^2/(4a)}\eta\left(q^{2a}\right)F_3(a,t;z),
\end{equation}
then\\
1)(Conjecture) The function $F_3(a,t;z)$ takes algebraic values, when $a,t\in\textbf{Q}^{*}_{+}$ and $z=r_1+i\sqrt{r_2}$, with $r_1$ rational and $r_2$ is positive rational.\\
2) If $a,t$ positive integers with $a>t$, then
$$
F_3(a,t;z)=Q^{\{3\}}_{\{a,t\}}(m(q))=
$$
\begin{equation}
=q^{-a/12+t^2/(4a)}\exp\left(-\sum^{\infty}_{n=1}q^n\sum_{\scriptsize
\begin{array}{cc}
AB=n\\
B\equiv \pm (a-t)(mod 2a)	
\end{array}\normalsize}\frac{(-1)^A}{A}\right).
\end{equation}
3) For the transformation of variables (220) holds
\begin{equation}
F_3(a',t';z')=\exp\left(\frac{-i\pi t^2 z}{2a}\right)F_3(a,t;z).
\end{equation}
\\
\textbf{Theorem 29.} $(Conjecture)$\\
When $a>t$ and $a,t$ positive rationals, the function $Q^{\{3\}}_{\{a,t\}}(x)$, takes algebraic numbers to algebraic numbers.\\
\\
\textbf{Notes.}\\
If $a,b,p$ are positive reals, with $a<a+b<p$ and $q=e(z)$, $Im(z)>0$, then the Ramanujan quantity $(RQ_4)$ is defined as (see [18]):
$$
RQ_4(a,b;p;z):=q^{(a^2-b^2)/(2p)-(a-b)/2}\frac{[a,p;q]_{\infty}}{[b,p,q]_{\infty}}=
$$
\begin{equation}
=q^{(a^2-b^2)/(2p)-(a-b)/2}\frac{(q^a;q^p)_{\infty}(q^{p-a};q^p)_{\infty}}{(q^b;q^p)_{\infty}(q^{p-b};q^p)_{\infty}}\textrm{, }|q|<1,
\end{equation}
where 
\begin{equation}
(a;q)_{\infty}:=\prod^{\infty}_{n=0}(1-aq^n)\textrm{ and }[a,p;q]_{\infty}:=\left(q^a;q^p\right)_{\infty} \left(q^{p-a};q^p\right)_{\infty}.
\end{equation}
Then using Jacobi triple product identity (relation (207) and related references), we get
\begin{equation}
RQ_4(a,b;p;z)=q^{(a^2-b^2)/(2p)-(a-b)/2}\frac{\vartheta_4\left(\left(\frac{p}{2}-a\right)\pi z;q^{p/2}\right)}{\vartheta_4\left(\left(\frac{p}{2}-b\right)\pi z;q^{p/2}\right)}.
\end{equation}
Here we consider also the case
\begin{equation}
RQ_3(a,b;p;z):=q^{(a^2-b^2)/(2p)-(a-b)/2}\frac{(-q^a;q^p)_{\infty}(-q^{p-a};q^p)_{\infty}}{(-q^p;q^p)_{\infty}(-q^{p-b};q^p)_{\infty}}.
\end{equation}
Then using Jacobi triple product identity we get
\begin{equation}
RQ_3(a,b;p;z)=q^{(a^2-b^2)/(2p)-(a-b)/2}\frac{\vartheta_3\left(\left(\frac{p}{2}-a\right)\pi z;q^{p/2}\right)}{\vartheta_3\left(\left(\frac{p}{2}-b\right)\pi z;q^{p/2}\right)}.
\end{equation}
From Theorem 28 relation (222) we have
$$
RQ_3(a,b;p;q)=
$$
$$
=q^{(a^2-b^2)/(2p)-(a-b)/2}\frac{q^{p/12-(p/2-a)^2/(2p)}}{q^{p/12-(p/2-b)^2/(2p)}}\frac{\eta(q^{p})}{\eta(q^p)}\frac{F_3\left(\frac{p}{2},\frac{p}{2}-a;z\right)}{F_3\left(\frac{p}{2},\frac{p}{2}-b;z\right)}=
$$
$$
=\frac{F_3\left(\frac{p}{2},\frac{p}{2}-a;z\right)}{F_3\left(\frac{p}{2},\frac{p}{2}-b;z\right)}.
$$
Hence we can state the next\\
\\
\textbf{Theorem 30.}\\
Suppose that $q=e(z)$, $Im(z)>0$ and $0<a<b<a+b\leq p$, where $a,b,p$ are reals. Then
\begin{equation}
RQ_4(a,b;p;q)=\frac{F_4\left(\frac{p}{2},\frac{p}{2}-a;z\right)}{F_4\left(\frac{p}{2},\frac{p}{2}-b;z\right)},
\end{equation}
where 
\begin{equation}
F_4\left(a,t;z\right):=\frac{\vartheta_4\left(\pi t z,e(az)\right)}{q^{a/12-t^2/(4a)}\eta\left(q^{2a}\right)}
\end{equation}
and
\begin{equation}
RQ_3(a,b;p;q)=\frac{F_3\left(\frac{p}{2},\frac{p}{2}-a;z\right)}{F_3\left(\frac{p}{2},\frac{p}{2}-b;z\right)}.
\end{equation}
Also for the transformation of variables 
\begin{equation}
p'=1/p\textrm{, }a'=2az/p\textrm{, }b'=2bz/p\textrm{, }z'=-1/(4z)\textrm{, }q'=e(-1/(4z)),
\end{equation}
we have the next modular identity
$$
RQ_3\left(p'-a',p'-b';2p';e\left(-\frac{1}{4z}\right)\right)=
$$
\begin{equation}
=\exp\left(-\pi i z \frac{a^2-b^2}{2p}\right)RQ_3\left(p-a,p-b;2p;e(z)\right).
\end{equation}
\\

Ramanujan has stated that (see [9] pg.21):\\
\\
\textbf{Theorem 31.}\\
Suppose that $q$, $a$ and $b$ are complex numbers with $\left|q\right|<1$, or that $q,a$, and $b$ are complex numbers with $a=bq^m$ for some integer $m$. Then
$$
U=U(a,b;q)=\frac{(-a;q)_{\infty}(b;q)_{\infty}-(a;q)_{\infty}(-b;q)_{\infty}}{(-a;q)_{\infty}(b;q)_{\infty}+(a;q)_{\infty}(-b;q)_{\infty}}=
$$  
\begin{equation}
=\frac{a-b}{1-q+}\frac{(a-bq)(aq-b)}{1-q^3+}\frac{q(a-bq^2)(aq^2-b)}{1-q^5+}\frac{q^2(a-bq^3)(aq^3-b)}{1-q^7+}\ldots.
\end{equation}
\\
\textbf{Theorem 32.}(see [19] Theorem A5 (Appendix))\\
If $0<a<p$ and $q=e(z)$, $Im(z)>0$, then holds the following continued fraction expansion
$$
\frac{\vartheta_3\left(\left(\frac{p}{2}-a\right)\pi z;q^{p/2}\right)}{\vartheta_4\left(\left(\frac{p}{2}-a\right)\pi z;q^{p/2}\right)}=\frac{F_3\left(\frac{p}{2},\frac{p}{2}-a;z\right)}{F_4\left(\frac{p}{2},\frac{p}{2}-a;z\right)}=
$$
\begin{equation}
=\frac{\sqrt{m^{*}_1\left(pz\right)}}{\textrm{dn}\left((p-2a)zK\left[pz\right],q^p\right)}
=-1+\frac{2}{1-U\left(q^a,-q^{p-a};q^p\right)}.
\end{equation}
Moreover if $a,p$ are positive integers, then
\begin{equation}
\log\left(-1+\frac{2}{1-U\left(q^a,-q^{p-a};q^p\right)}\right)=2\sum^{\infty}_{n=1}q^n\sum_{\scriptsize
\begin{array}{cc}
AB=n\\
A\equiv1(2)\\
B\equiv \pm a(p)	
\end{array}}\frac{1}{A}.
\end{equation}
Here we have use
\begin{equation}
K(w)=\frac{\pi}{2}\cdot  {}_2F_1\left(\frac{1}{2},\frac{1}{2};1;w^2\right)\textrm{, }|w|<1
\end{equation}
and
\begin{equation}
K[z]:=K(m(q))\textrm{, }q=e(z)\textrm{, }Im(z)>0,
\end{equation}
where
\begin{equation}
m^{*}(z)=m(q)=\left(\frac{\vartheta_2(0,q^{1/2})}{\vartheta_3(0,q^{1/2})}\right)^2\textrm{, }m_1^{*}(z)=\sqrt{1-m^{*}(z)^2}.
\end{equation}
\\
\textbf{Theorem 33.}\\
If $q=e(z)$, $Im(z)>0$ and $a,b,p$ reals such that $a<b<a+b\leq p$, then
\begin{equation}
\frac{RQ_3(a,b;p;q)}{RQ_4(a,b;p;q)}=\left(\frac{\textrm{dn}\left((p-2a)zK\left[pz\right],q^p\right)}{\textrm{dn}\left((p-2b)zK\left[pz\right],q^p\right)}\right)^{-1}.
\end{equation}
\\
\textbf{Theorem 33.1}\\
Assume that $p,a$ are  positive integers. Then
$$
\theta(a,p;z)=q^{p/8+a^2/(2p)-a/2}\sum^{\infty}_{n=-\infty}(-1)^nq^{pn^2/2+(p/2-a)n}\textrm{, }q=e^{i\pi z}\eqno{(241.1)}
$$
is a modular form of weight $1/2$ in $\Gamma(2p)$. That is if $a_1,b_1,c_1,d_1$ are integers such that $a_1,d_1\equiv1(\textrm{mod}2p)$, $b_1,c_1\equiv0(\textrm{mod}2p)$ and $a_1d_1-b_1c_1=1$, we get 
$$
\theta\left(\frac{a_1z+b_1}{c_1z+d_1}\right)=\epsilon\sqrt{c_1z+d_1}\theta(z)\textrm{, }Im(z)>0,\eqno{(241.2)}
$$ 
where $\epsilon$ depends only on $a_1,b_1,c_1,d_1$ and  $\epsilon^{24}=1$.\\ 
\\

For example set $a=16$, $p=128$, then
$$
\phi(z)=q^9\sum^{\infty}_{n=-\infty}(-1)^nq^{64n^2+48n}=\sum^{\infty}_{n=-\infty}(-1)^nq^{(8n+3)^2}=\sum^{\infty}_{n=-\infty}\chi_0(n)q^{n^2}=
$$
$$
=\sum_{n\equiv 3(\textrm{\scriptsize mod\normalsize})8}(-1)^{\frac{n-3}{8}}q^{n^2}=
$$
$$
=q^9-q^{25}-q^{121}+q^{169}+q^{361}-q^{441}-q^{729}+q^{841}+\ldots
$$
is a modular form in $\Gamma_1(128)$ of weight $1/2$ and all coefficients are non zero only at  $3(\textrm{mod}8)$. Also $\chi_0(n)=(-1)^{\frac{n-3}{8}}$. More general
$$
\theta(a,p;8pz)=\sum^{\infty}_{n=-\infty}(-1)^nq^{(2np+p-2a)^2}=\sum_{\scriptsize \begin{array}{cc}
n\in\textbf{\scriptsize Z}\\
n\equiv p-2a(\textrm{\scriptsize mod\normalsize}2p)
\end{array}
\normalsize}(-1)^{\frac{n-(p-2a)}{2p}}q^{n^2}.
$$
Hence
$$
\theta(a,p;z)=\sum_{\scriptsize \begin{array}{cc}
n\in\textbf{\scriptsize Z}\\
n\equiv p-2a(\textrm{\scriptsize mod\normalsize}2p)
\end{array}\normalsize
}(-1)^{\frac{n-(p-2a)}{2p}}q^{n^2/(8p)}.\eqno{(241.3)}
$$
If we assume that $X(n,m)$ is any bouble arithmetical function and set 
$$
\sum^{\infty}_{n,m=-\infty}X(n,m)q^{n^2+m^2}=\sum^{\infty}_{n=0}R(n)q^{n},\eqno{(241.4)}
$$
$$
\textrm{Sym}^{*}X\left(n,m\right):=\frac{1}{2}\left(X(n,m)+X(m,n)+X(-n,-m)+X(-m,-n)\right),\eqno{(241.5)}
$$
and
$$
A(k,n)=\left\{\begin{array}{cc}
\frac{1}{2}\textrm{Sym}^*X\left(-\frac{\sqrt{2n}}{2},\frac{\sqrt{2n}}{2}\right)\textrm{, if }\sqrt{2n}\in\textbf{N}\textrm{ and }k=\sqrt{2n}\\
\textrm{Sym}^*X\left(x^{-}_{k,n},x^{+}_{k,n}\right)\textrm{, if }\sqrt{2n}\textrm{ not in }\textbf{N}
\end{array}\right\},\eqno{(241.6)}
$$
where 
$$
x^{\pm}_{k,n}=\frac{1}{2}\left(\pm k-\sqrt{2n-k^2}\right).
$$
Then 
$$
R(n)=\sum_{\scriptsize\begin{array}{cc} 
0\leq |k|\leq\left[\sqrt{2n}\right]\\
2n-k^2=l^2\geq0
\end{array}\normalsize}A(k,n).
$$
Hence\\
\\
\textbf{Theorem 33.2}\\
In general holds
$$
\sum^{\infty}_{n,m=-\infty}X(n,m)q^{n^2+m^2}=\sum^{\infty}_{n=0}q^n\sum_{\scriptsize\begin{array}{cc} 
0\leq |k|\leq\left[\sqrt{2n}\right]\\
2n-k^2=l^2\geq0
\end{array}\normalsize}A(k,n).\eqno{(241.7)}
$$
\\

Setting in (241.7) 
$$
X(n,m)=\chi_{p-2a,2p}(n)(-1)^{\frac{n-(p-2a)}{2p}}\chi_{p-2a,2a}(m)(-1)^{\frac{m-(p-2a)}{2p}},
$$
where $\chi_{a,b}(n)=1$ if $n$ is integer of the form $n\equiv a(\textrm{mod}b)$, and 0 else, we get
$$
A(k,n)=A_0(k,n)=
$$
$$
-i^{\frac{4a-\sqrt{2n-k^2}}{p}}\chi_{p-2a,2p}\left(\frac{-k-\sqrt{2n-k^2}}{2}\right)\chi_{p-2a,2p}\left(\frac{k-\sqrt{2n-k^2}}{2}\right)-
$$
$$
-i^{\frac{4a+\sqrt{2n-k^2}}{p}}\chi_{p-2a,2p}\left(\frac{-k+\sqrt{2n-k^2}}{2}\right)\chi_{p-2a,2p}\left(\frac{k+\sqrt{2n-k^2}}{2}\right).
$$
Hence we have 
$$
\left(\sum^{\infty}_{n=-\infty}(-1)^nq^{pn^2/2+(p/2-a)n}\right)^2=
$$
$$
=q^{-a^2/p-p/4+a}\sum^{\infty}_{n=0}q^{n/(8p)}\sum_{\scriptsize\begin{array}{cc} 
0\leq |k|\leq\left[\sqrt{2n}\right]\\
2n-k^2=l^2\geq0
\end{array}\normalsize}A_0(k,n),
$$
where 
$$
A_0(k,n)=-i^{\frac{4a-l}{p}}\chi_{p-2a,2p}\left(-\frac{k+l}{2}\right)\cdot\chi_{p-2a,2p}\left(\frac{k-l}{2}\right)-
$$
$$
-i^{\frac{4a+l}{p}}\chi_{p-2a,2p}\left(-\frac{k-l}{2}\right)\cdot\chi_{p-2a,2p}\left(\frac{k+l}{2}\right)
$$
and $l=\sqrt{2n-k^2}$. Hence we get the next:\\
\\
\textbf{Theorem 33.3}\\
If $a,p$ are integers with $p>0$ and $p>2|a|$, then
$$
\left(\sum^{\infty}_{n=-\infty}(-1)^nq^{pn^2/2+(p/2-a)n}\right)^2
=q^{-a^2/p-p/4+a}\sum^{\infty}_{n=0}C(a,p,n)q^{n/(8p)},\eqno{(241.8)}
$$
where 
$$
C(a,p,n)=-\sum_{\scriptsize\begin{array}{cc} 
0\leq |k|\leq\left[\sqrt{2n}\right]\\
2n-k^2=l^2\geq1\\
k\equiv 0(\textrm{mod}2p)\\
l-k\equiv(2p\pm 4a)(\textrm{mod}4p)
\end{array}\normalsize}(-1)^{\frac{4a\mp l}{2p}}.\eqno{(241.9)}
$$
\\
\textbf{Theorem 33.4}\\
If $\chi(n)$ is full multiplicative function in $\textbf{Z}$ with $\chi(0)=0$, $\chi(1)=1$, then
$$
\left(\sum^{\infty}_{n=-\infty}\chi(n)q^{n^2}\right)^2=\sum^{\infty}_{n=0}C_{\chi}(n)q^{n/16},\eqno{(241.10)}
$$
where 
$$
C_{\chi}(n)=\sum_{\scriptsize\begin{array}{cc} 
0\leq |m|\leq\left[\sqrt{2n}\right]\\
2n-m^2=0\\
m\equiv 0(\textrm{mod}8)
\end{array}\normalsize}\chi\left(\frac{-m^2}{64}\right)+2\sum_{\scriptsize\begin{array}{cc} 
0\leq |m|\leq\left[\sqrt{2n}\right]\\
2n-m^2=l^2\geq1\\
m\equiv 0(\textrm{mod}4)\\
l-m\equiv 0(\textrm{mod}8)
\end{array}\normalsize}\chi\left(\frac{l^2-m^2}{64}\right).\eqno{(241.11)}
$$

\section{Further conjectures and notes}

In this, and the next paragraph, we consider notes and conjectures, mixed together for further study. The reader must read it carefuly and decide what will keep and what will left behind, althought the results have been checked numericaly and some of them have been proved. I want also to mention that this research it is not about, if the conjecture of paragraph 2 is correct or not. It is in one way, to study general properties of the $A(a,p,q)$ function, and in other way, to study such general phenomena in $q-$series. That is why we consider general transformations such in [15], and paragraphs 4,5,6,7 of the present article. What will hapen if we break the law of ``quadraticity'' and go to higher forms? My search is to try address this problem. This will continue in future works (see preperations [23],[24]).\\
\\   

Assume now the notation of (238),(239),(240). In [8] we have considered for $|q|<1$ the function (here we make some modifications in $\tau^{*}$)
\begin{equation}
\tau^{*}(a,p;z):=\tau(a,p;q):=q^{C}\frac{\left[a,p;q^2\right]_{\infty}}{\left[a,p;q\right]_{\infty}}\textrm{, }C=\frac{a^2}{2p}-\frac{a}{2}+\frac{p}{12},
\end{equation}
with
\begin{equation}
\tau^{*}(a,p;z)=\tau^{*}(np\pm a,p;z)\textrm{, }\forall n\in\textbf{Z}.
\end{equation}
Hence\\ 
\\
\textbf{Conjecture 1.}\\
We assume that exists holomorphic function $\tau_e^{*}(a,p,z)=\tau_e\left(q_1,p,q\right)$, where $q_1=e(a/p)$, $q=e(z)$, $Im(a)>0$, $Im(z)>0$ and $0<p-$parameter, such that
\begin{equation}
\tau^{*}(a,p;z)=\tau_e\left(e(a/p),p;e(z)\right)+\tau_e\left(e(-a/p),p;e(z)\right). 
\end{equation}
\\

Hence if we assume Conjecture 1, i.e. that $\tau^{*}(a,p;z)$ is $p-$periodic and ''even'' with respect to $a$, then we get the validity of (243). Also 
\begin{equation}
\left(\frac{\partial\tau^{*}(a,p;z)}{\partial a}\right)_{a\in\frac{p}{2}\textbf{\scriptsize Z\normalsize}}=0.
\end{equation}
Hence if 
\begin{equation}
\tau_0(a;q):=\tau(a,1,q)=q^{1/12-a/2+a^2/2}\frac{[2a,2;q]_{\infty}}{[a,1;q]_{\infty}},
\end{equation}
then
\begin{equation}
\tau_0(n\pm a;q)=\tau_0(a;q)\textrm{, }\forall n\in\textbf{Z}
\end{equation}
and
\begin{equation}
\left(\frac{\partial \tau_0(a;q)}{\partial a}\right)_{a\in\frac{1}{2}\scriptsize\textbf{Z}\normalsize}=0
\end{equation}
and
\begin{equation}
\left(\frac{\partial \tau_0(a;q)}{\partial a}\right)_{a=n+t}=\left(\frac{\partial \tau_0(a;q)}{\partial a}\right)_{a=t}\textrm{, }\forall n\in\textbf{Z}\textrm{, }t\in\textbf{R}.
\end{equation}
Moreover from [19] Appendix, we have 
\begin{equation}
\frac{\left[a,p;q^2\right]_{\infty}}{\left[a,p;q\right]^2_{\infty}}=\frac{\vartheta_3\left((p/2-a)\pi z;q^{p/2}\right)}{\vartheta_4\left((p/2-a)\pi z;q^{p/2}\right)}=\frac{\textrm{dn}\left((p-2a)zK[pz],q^{p}\right)}{\sqrt{m_1^{*}\left(pz\right)}}.
\end{equation} 
Considering the above foundings we get the next evaluation\\
\\
\textbf{Theorem 34.}\\
If $q=e(z)$, $Im(z)>0$ and $a,p>0$, then with the notation of (238), (239), (240), (242), we have
\begin{equation}
q^{\frac{p}{12}-\frac{a}{2}+\frac{a^2}{2p}}\left[a,p;q\right]_{\infty}=A(a,p;q)=\frac{\sqrt{m_1^{*}(pz)}\tau^{*}(a,p;z)}{\textrm{dn}\left((p-2a)zK[pz],q^p\right)}.
\end{equation}
\\
And\\
\\
\textbf{Theorem 35.}(Conjercture)\\
If $n\in\textbf{Z}$ and $q=e(z)$, $Im(z)>0$, then
\begin{equation}
A(np\pm a,p;q)=\pm (-1)^nA(a,p,q)
\end{equation}
and
\begin{equation}
\left(\frac{\partial A(a,p;q)}{\partial a}\right)_{a=(n+1/2)p}=0
\end{equation}
and
\begin{equation}
\left(\frac{\partial A(a,p;q)}{\partial a}\right)_{a=np}=-(-1)^n2\pi iz q^{p/24}\eta_D(pz).
\end{equation}
\\
\textbf{Proof.}\\
From Theorem 34 and relations (243),(245) and 
$$
\frac{\textrm{dn}\left(((2n-1)p+2a)zK[pz],q^p\right)}{\textrm{dn}\left((p-2a)zK[pz],q^p\right)}=(-1)^n\textrm{, }n\in\textbf{Z},
$$
we get the three results.\\
\\
\textbf{Theorem 36.}\\
Assume $q=e^{\pi i z}$, $Im(z)>0$ and set 
\begin{equation}
\tau_0^{*}(a):=\tau^{*}_0(a;z)=\tau_0(a;q)=q^{1/8-a/2+a^2/2}\prod^{\infty}_{n=0}\left(1+q^{n+a}\right)\left(1+q^{n+1-a}\right)
\end{equation}
and
\begin{equation}
\tau_1^{*}(a):=\tau^{*}_1(a;z)=\tau_1(a;q)=q^{1/8-a/2+a^2/2}\prod^{\infty}_{n=0}\left(1-q^{n+a}\right)\left(1-q^{n+1-a}\right).
\end{equation}
Then\\
1) 
\begin{equation}
\tau^{*}_0(a+1)=\tau^{*}_0(a)
\end{equation} 
and 
\begin{equation}
\tau^{*}_0\left(a+\frac{2}{z}\right)=-e^{2\pi i a}e^{2\pi i/z}\tau^{*}_0(a).
\end{equation}
Also\\
2)
\begin{equation}
\tau^{*}_1(a+2)=\tau^{*}_1(a)
\end{equation} 
and 
\begin{equation}
\tau^{*}_1\left(a+\frac{2}{z}\right)=-e^{2\pi i a}e^{2\pi i/z}\tau^{*}_1(a).
\end{equation}
3) The function
\begin{equation}
\tau_{10}^{*}(a):=\frac{\tau^{*}_1(a)}{\tau^{*}_0(a)},
\end{equation}
is double periodic, with periods $2$ and $2/z$ and hence it is an elliptic function.\\
4) 
\begin{equation}
\tau^{*}_1(a;2z)=\tau^{*}_1(a;z)\tau^{*}_0(a;z).
\end{equation}
5)
\begin{equation}
\tau_0^{*}(a)=\frac{q^{1/8-a/2+a^2/2}}{\eta_D\left(z/2\right)}\sum^{\infty}_{n=-\infty}q^{n^2/2+(1/2-a)n}\textrm{, }q=e^{\pi i z}.
\end{equation}
\begin{equation}
\tau_1^{*}(a)=\frac{q^{1/8-a/2+a^2/2}}{\eta_D\left(z/2\right)}\sum^{\infty}_{n=-\infty}(-1)^nq^{n^2/2+(1/2-a)n}\textrm{, }q=e^{\pi i z}.
\end{equation}
6)
\begin{equation}
A(a+2,1;z)=A(a,1;z)
\end{equation}
and
\begin{equation}
A\left(a+\frac{2}{z},1;z\right)=-e^{2\pi i a}e^{2\pi i/z}A\left(a,1;z\right).
\end{equation}
\\
\textbf{Proof.}\\
See [2].\\
\\
\textbf{Theorem 37.}\\
If $q=e(z)$, $Im(z)>0$ and
\begin{equation}
\tau^{*}(a,p;z)=q^{p/12-a/2+a^2/(2p)}\left(-q^a;q^p\right)_{\infty}\left(-q^{p-a};q^p\right)_{\infty},
\end{equation}
then for $a$ complex and $2p\in\textbf{N}$, we have
\begin{equation}
\tau^{*}\left(a+p,p;z\right)=\tau^{*}(a,p;z)
\end{equation}
and
\begin{equation}
\tau^{*}\left(a+\frac{2p}{z},p;z\right)=e^{2\pi i (p+2a)}e^{4\pi ip/z}\tau^{*}(a,p;z).
\end{equation}
Also
\begin{equation}
A(a+2p,p;z)=A(a,p;z)
\end{equation}
and
\begin{equation}
A\left(a+\frac{2p}{z},p;z\right)=e^{2\pi i (p+2a)}e^{4\pi i p/z}A(a,p;z).
\end{equation}
\\
\textbf{Proof.}\\
See [2].\\
\\
\textbf{Corollary.}\\
We define
\begin{equation}
J(a,p;z):=\frac{\sqrt{m_1(pz)}}{\textrm{dn}\left((p-2a)z K[pz],q^p\right)}.
\end{equation}
If $a\in\textbf{C}$, $2p\in\textbf{N}$ and $Im(z)>0$, then
\begin{equation}
J(a+2p,p;z)=J(a,p;z)
\end{equation}
and
\begin{equation}
J\left(a+\frac{2p}{z},p;z\right)=J(a,p;z).
\end{equation}
\\
\textbf{Proof.}\\
Easy.\\
\\

Assume again $a\in\textbf{C}$, $2p\in\textbf{N}$ and $q=e^{2\pi i z}$, $Im(z)>0$. Also again define
\begin{equation}
A_0(a,p;z):=A^{*}(a,p;q)=q^{p/12-a/2+a^2/(2p)}\left(-q^a;q^p\right)_{\infty}\left(-q^{p-a};q^p\right)_{\infty}.
\end{equation}
Then the function defined as
\begin{equation}
A_1(a;z):=A_0\left(a,1;-1/z\right)=\tau^{*}_{0}(a,-2/z),
\end{equation}
satisfy the relations
\begin{equation}
A_1(a+1;z)=A_1(a;z)
\end{equation}
and
\begin{equation}
A_1\left(a+z;z\right)=e^{2\pi i (1/2-a-z/2)}A_1(a;z).
\end{equation}
Hence the function $A_1(a;z)^2$ is nearly a Jacobi form (see [22]). This follows from 
\begin{equation}
A_1(a+m;z)=A_1(a;z)\textrm{, }m\in\textbf{Z}
\end{equation} 
and 
\begin{equation}
A_1(a+lz;z)=e^{2\pi i l(1/2-a-l z/2)}A_1(a;z)\textrm{, }l\in\textbf{Z}.  
\end{equation}
Hence we have
\begin{equation}
A_1(a+lz+m;z)^2=e^{-2 \pi i(l^2 z+2la)}A_1(a;z)^2,
\end{equation}
Also if we set 
\begin{equation}
\tau_{00}(a,z):=\tau^*_0\left(a,-2/z\right),
\end{equation} 
then from 
\begin{equation}
\frac{\vartheta_3\left(\pi t z;e\left(pz/2\right)\right)}{q^{p/24-t^2/(2p)}\eta\left(q^p\right)}=A^{*}\left(a,p;q\right)\textrm{, }t=p/2-a,
\end{equation}
and the modular relations (219),(220) one can arrive to
\begin{equation}
\tau_{00}\left(a;\frac{a_0z+b_0}{c_0z+d_0}\right)^{12}=\tau_{00}\left(a;z\right)^{12}\textrm{, when }a\in\textbf{Z}\textrm{, }Im(z)>0
\end{equation}
and $a_0d_0-b_0c_0=1$, $c_0,b_0\equiv0(\textrm{mod}2)$ and $a_0,d_0\equiv1(\textrm{mod}2)$. Hence  $\tau_{00}(a;z)^{6}$ is a weight 0 modular form on $\Gamma(2)$. Hence we have the next\\
\\
\textbf{Theorem 38.}(Conjecture)\\
When $a$ is integer, the function $\tau_{00}(a;z)^{12}$, is a weight $0$ modular form on $\Gamma(2)$ (as function of  $z$). Also when $a\in\textbf{C}$ and $Im(z)>0$, then $A_1(a;z)=A_0(a,1;-1/z)=\tau_{00}(a;z/2)$ have the properties
\begin{equation}
A_{1}(a+lz+m;z)^{12}=e^{-2 \pi i 6(l^2 z+2la)}A_{1}(a;z)^{12}\textrm{, }\forall l,m\in\textbf{Z}.
\end{equation}
and
\begin{equation}
A_{1}\left(\frac{a}{c_0z+d_0};\frac{a_0z+b_0}{c_0z+d_0}\right)^{12}=\exp\left(\frac{2\pi i 6c_0 a^2}{c_0z+d_0}\right)A_{1}\left(a;z\right)^{12},
\end{equation}
when $a_0,d_0,b_0,c_0\in\textbf{Z}$, $a_0d_0-b_0c_0=1$, $c_0,b_0\equiv0(\textrm{mod}2)$ and $a_0,d_0\equiv1(\textrm{mod}2)$. Hence $A_1(a;z)^{12}$ is a Jacobi modular form of weight $0$ and level $N=6$ on $\Gamma(2)$.\\ 
\\

Now $A_0(a,1;z)=\tau_0^{*}(a;2z)$ and from Theorem 36 relation (263), setting where $q\rightarrow q^p$, we get ($q=e^{2\pi i z}=e(z)$):
$$
\tau_0^{*}(a;2zp)=\frac{q^{p/8-ap/2+a^2p/2}}{\eta_D\left(zp\right)}\sum^{\infty}_{n=-\infty}q^{pn^2/2+(p/2-pa)n}.
$$   
Setting then $a\rightarrow a/p$, we get
$$
\tau_0^{*}(a/p;2zp)=A(a,p;q)=\tau^{*}(a,p;z)=\frac{q^{p/8-a/2+a^2/(2p)}}{\eta_D\left(zp\right)}\sum^{\infty}_{n=-\infty}q^{pn^2/2+(p/2-a)n}.
$$ 
Hence we get the next\\
\\
\textbf{Theorem 39.}(Conjecture)\\
Assume $q=e(z)$, $Im(z)>0$. Then if $p>0$ and $a$ complex number, we have
\begin{equation}
\sum^{\infty}_{n=-\infty}q^{pn^2/2+(p/2-a)n}=q^{-p/8+a/2-a^2/(2p)}\eta_{D}(zp)\tau_0^*(a/p;2zp), 
\end{equation} 
where $\tau_0^{*}(a;z)=A_0(a,1;z/2)$. Also, when $a\in\textbf{Z}$, the function  $\tau_0^*\left(a;z\right)^{12}$ is a weight 0 modular form on $\Gamma(2)$.\\
\\

But it is known that if $J_0(z)$, $Im(z)>0$, denotes the Klein invariant ($J_0(0)=1$), then $J_0(2z)$ is also a weight 0 modular form but now in $\Gamma(2)$. Hence when $a$ is integer $\tau_0^*\left(a;z\right)^{12}$ is an algebraic function of $J_0(2z)$. By this way we can write $\tau_0^*\left(a;z\right)^{12}=f\left(J_0(2z)\right)
$. Inverting this last equation we get 
$$
\tau^{*}_{0}\left(a;\frac{1}{2}J^{(-1)}_0(z)\right)^{12}=f(z)
$$ 
and when $J_0(\frac{1}{2}\rho)=t$ is positive rational, the function $f(t)$ is always solution of a sextic equation with integer coefficients. Hence we get the next:\\
\\
\textbf{Theorem 40.} (Conjecture)\\
If $a$ is poisitive rational $a_1/a_2$, with $a_1,a_2\in\textbf{N}$, $a_2>1$ and $(a_1,a_2)=1$, then  $\tau^{*}_0\left(a_1/a_2;z\right)^{4P(a_2)}$ is a modular form of weight 0 on $\Gamma\left(2P(a_2)\right)$ and exists algebraic function $f(x)$ (depending on $a_2$), such that  
\begin{equation}
\tau^{*}_0\left(\frac{a_1}{a_2};z\right)^{4P(a_2)}=f\left(J_0\left(2P(a_2)z\right)\right)\textrm{, }Im(z)>0.
\end{equation}
In the case $a_2=1$, assume $J_0(2\rho)=t$ is positive rational, then the function $U=\tau^*_0\left(a;\rho\right)^{12}$ is solution of a sextic equation with rational coefficients. The function $P(n)$ is defined for integer $n>0$ as: 
\begin{equation}
P(n)=\prod_{\scriptsize
\begin{array}{cc}
1<p|n\\
p-prime
\end{array}\normalsize}p.
\end{equation}
Also $J_0(z)$ denotes the Klein's $J_0-$invariant.\\
\\
\textbf{Theorem 41.}\\
Assume that $p,a$ are positive integers. Then
\begin{equation}
\theta(a,p;z)=q^{p/8+a^2/(2p)-a/2}\sum^{\infty}_{n=-\infty}q^{pn^2/2+(p/2-a)n}\textrm{, }q=e^{\pi i z},
\end{equation}
is a modular form of weight $1/2$ in $\Gamma(2p)$. That is if $a_1,b_1,c_1,d_1$ are integers such that $a_1,d_1\equiv 1(\textrm{mod }2p)$, $b_1,c_1\equiv0 (\textrm{mod }2p)$ and $a_1d_1-b_1c_1=1$, we get
\begin{equation}
\theta\left(\frac{a_1 z+b_1}{c_1 z+d_1}\right)=\epsilon \sqrt{c_1 z+d_1}\theta(z)\textrm{, }Im(z)>0,
\end{equation}
where $\epsilon$ depends only on $a_1,b_1,c_1,d_1$ and $\epsilon^{24}=1$.

\section{Table for $\theta_4(p/2,p/2-a;q)$}

Here we give a table of evaluations, which does not include Theorems 1 and 2, for certain lower values $a,p$ of theta functions.\\
\\
\textbf{1.}
\begin{equation}
\sum^{\infty}_{n=-\infty}(-1)^nq^{3/2n^2+n/2}=q^{1/12}\eta(q^3)Q(k_r^2)
\end{equation}
The polynomial which relates $u=A(1,3,q)^{12}$ with $v=m(q)=k_r^2$ is
$$
u^4 v^5-4 u^4 v^4+6 u^4 v^3-4 u^4 v^2+u^4 v-16 u^3 v^6+84 u^3 v^5
-12480 u^3 v^4-
$$
$$
-40712 u^3 v^3-12480 u^3 v^2+84 u^3 v-16 u^3+196830 u^2 v^5-787320 u^2 v^4+1180980 u^2 v^3-
$$
$$
-787320 u^2 v^2+196830 u^2 v+19131876 u v^5-76527504 u v^4+114791256 u v^3-
$$
$$
-76527504 u v^2+19131876 u v+387420489 v^5-1549681956 v^4+2324522934 v^3-
$$
\begin{equation}
-1549681956 v^2+387420489 v=0
\end{equation}
\\
\textbf{2.}
\begin{equation}
\sum^{\infty}_{n=-\infty}(-1)^nq^{3n^2-5n}=q^{-11/6}\eta(q^6)Q(k_r)
\end{equation} 
The polynomial equation which relates $u=A(8,6;q)^6$ with $v=m(q)^{1/2}=k_r$ is 
$$
u^8 v^4-u^8 v^2+16 u^6 v^6-24 u^6 v^4-24 u^6 v^2+16 u^6-486 u^4 v^4+
$$
\begin{equation}
+486 u^4 v^2-19683 v^4+19683 v^2=0
\end{equation} 
\\
\textbf{3.}
\begin{equation}
\sum^{\infty}_{n=-\infty}(-1)^nq^{3n^2+4n}=q^{-13/12}\eta(q^6)Q(k_r)
\end{equation} 
The polynomial equation which relates $u=A(-1,6,q)^6$, with $v=m(q)^{1/2}=k_r$ is
\begin{equation}
u^4 v^3-u^4 v+16 u^3 v^2-18 u^2 v^3+18 u^2 v+4 u v^4-8 u v^2+4 u+v^3-v=0
\end{equation}
\\
\textbf{4.}
\begin{equation}
\sum^{\infty}_{n=-\infty}(-1)^nq^{4n^2+6n}=q^{-23/12}\eta(q^8)Q(k_r)
\end{equation}
The polynomial equation which relates $u=A(-2,8;q)^{12}$, with $v=m(q^2)^2=k_{4r}^4$ is
\begin{equation}
-u^4v-64 u^2 v+256 v^2-512 v+256=0
\end{equation}
\\
\textbf{5.}
\begin{equation}
\sum^{\infty}_{n=-\infty}(-1)^nq^{5/2n^2+3/2n}=q^{-1/60}\eta(q^5)Q(\eta_5(q^4)^5)  
\end{equation}
where 
\begin{equation}
\eta_5(q)=\frac{1}{2}\left(-1-h_5(q)+\sqrt{5+2h_5(q)+h_5(q)^2}\right)
\end{equation}
and 
\begin{equation}
h_5(q)=\frac{\eta(q^{1/5})}{q^{1/5}\eta(q^5)}
\end{equation} 
The polynomial equation which relates $u=A(1,5,q^2)^{15}$ with $v=\eta_5(q^4)^5$ is
$$
u^4+v^{11}+55 v^{10}+1205 v^9+13090 v^8+69585 v^7+134761 v^6-69585 v^5+
$$
\begin{equation}
+13090 v^4-1205 v^3+55 v^2-v=0
\end{equation}
Note that
$$
h_5(q^5)=\frac{\eta(q)}{q\eta(q^{25})}=\frac{1}{\sqrt{M_5(q)M_5(q^5)}}\left(\frac{m(q)}{m(q^{25})}\right)^{1/24}\left(\frac{m^{*}(q)}{m^{*}(q^{25})}\right)^{1/6}
$$
is function of $k_r=\sqrt{m(q)}$, since (see [14])
\begin{equation}
\left(5M_5(q)-1\right)^5\left(1-M_5(q)\right)=256\cdot m(q)m^{*}(q)M_5(q)
\end{equation}
and $m^{*}(q)=1-m(q)$.\\
\\
For tables of singular modulus one can see [10],[11],[14].

\[
\]

\centerline{\textbf{References}}

[1]: M. Abramowitz and I.A. Stegun. ''Handbook of Mathematical Functions''. 9th ed., New York: Dover Pub. (1972)
\\

[2]: J.V. Armitage and W.F. Eberlein. ''Elliptic Functions''. 2006 ed., New York: Cambridge University, (2006) 
\\

[3]: N.D. Bagis. ''On a general sextic equation solved by the Rogers-Ramanujan continued fraction''. arXiv:1111.6023v2, (2012)
\\

[4]: N.D. Bagis. ''Evaluation of Fifth Degree Elliptic Singular Moduli''. arXiv:1202.6246v3, (2015)
\\

[5]: Toshijune Miyake. ''Modular Forms''. Springer Verlag, Berlin, Heidelberg, (1989).
\\

[6]: N.D. Bagis. ''Parametric Evaluations of the Rogers Ramanujan Continued Fraction''. International Journal of Mathematics and Mathematical Sciences. Vol(2011). doi:101155/2011/940839.
\\

[7]: N.D. Bagis and M.L. Glasser. ''Conjectures on the evaluation of alternative modular bases and formulas approximating 1/$\pi$''. Journal of Number Theory, Vol(132), 2012 pp. 2353-2370. 
\\

[8]: N.D. Bagis and M.L. Glasser. ''Jacobian Elliptic Functions, Continued Fractions and Ramanujan Quantities''. arXiv:1001.2660v1, (2010)
\\ 

[9]: B.C. Berndt. ''Ramanujan's Notedbooks Part III''. 1991 ed., New York: Springer-Verlag. (1991) 
\\

[10]: B.C. Berndt. ''Ramanujan's Notedbooks Part V''. 1998 ed., New York: Springer-Verlag. (1998) .
\\

[11]: J.M. Borwein and P.B. Borwein. ''Pi and the AGM''. 1987 ed., New York: John Wiley and Sons Inc., (1987). 
\\

[12]: D. Broadhurst. ''Solutions by Radicals at Singular Values $k_N$ from New Class Invariants for $N\equiv3mod8$''. arXiv:0807.2976 (math-phy), (2008). 
\\

[13]: E.T. Whittaker and G.N. Watson. ''A course on Modern Analysis''. 4th ed., Cambridge University Press. (1927).
\\

[14]: J.M. Borwein, M.L. Glasser, R.C. McPhedran, J.G. Wan, I.J. Zucker. ''Lattice Sums Then and Now''. Cambridge University Press. New York, (2013).
\\

[15]: N.D. Bagis. ''On Generalized Integrals and Ramanujan-Jacobi Special Functions''. arXiv:1309.7247v3, (2015).
\\

[16]: G.E. Andrews. ''Number Theory''. Dover Publications, New York. (1994).
\\

[17]: T. Apostol. ''Introduction to Analytic Number Theory''. Springer Verlang, New York, Berlin, Heidelberg, Tokyo. (1974).
\\

[18]: N.D. Bagis. ''Generalizations of Ramanujan Continued Fractions''. arXiv:1107.2393v2 [math.GM], (2012).
\\

[19]: N.D. Bagis. ''Evaluations of Series Related to Jacobi Elliptic Functions''. arXiv:1803.09445v2 [math.NT], (2019).
\\

[20]: B.C. Berndt. ''Ramanujan`s Notebooks Part II''. Springer-Verlag, New York. (1989).
\\

[21]: N.D. Bagis. ''Evaluations of certain theta functions in Ramanujan theory of alternative modular bases''. arXiv:1511.03716v2 [math.GM], (2017).
\\

[22]: A. Dadholkar, S. Murthy, D. Zagier. ''Quantum Black Holes, Wall Crossing and Mock Modular Forms''. arXiv:1208.4074v2 [hep-th], (2014).
\\

[23]: N.D. Bagis. ''$q-$series Related with Higher Forms''. arXiv:2006.16005v4 [math.GM], (2021).
\\

[24]: N.D. Bagis. ''On the numbers that are sums of three cubes''. arXiv:2009.11972v1 [math.GM], (2020).

\end{document}